\pdfoutput=1
\documentclass[a4paper]{article}
    \usepackage[a4paper,left=1in,right=1in,top=1in,bottom=1in]{geometry}

    \usepackage{titlesec}
    \usepackage[title]{appendix}
    
    \titleformat{\section}[frame]
    {\normalfont}
    {\filright
    \footnotesize
    \enspace SECTION \thesection\enspace}
    {8pt}
    {\Large\bfseries\filcenter}

    \usepackage{amsmath}
    \usepackage{amsthm}
    \usepackage{amssymb} 

    \usepackage{tikz-cd,tikz}
    \usetikzlibrary{arrows,backgrounds}

    \usepackage{mathtools}
    \usepackage{stmaryrd}

    \usepackage[utf8]{inputenc}
    \usepackage[T1]{fontenc}
    \usepackage{libertine} 
    \usepackage[libertine]{newtxmath}

    \usepackage[style=numeric,backref=true,giveninits=true,doi=false,url=false,isbn=false]{biblatex}
    
    \newbibmacro{string+doiurlisbn}[1]{%
      \iffieldundef{doi}{%
        \iffieldundef{url}{%
          \iffieldundef{isbn}{%
            \iffieldundef{issn}{%
              #1%
            }{%
              \href{http://books.google.com/books?vid=ISSN\thefield{issn}}{#1}%
            }%
          }{%
            \href{http://books.google.com/books?vid=ISBN\thefield{isbn}}{#1}%
          }%
        }{%
          \href{\thefield{url}}{#1}%
        }%
      }{%
        \href{http://dx.doi.org/\thefield{doi}}{#1}%
      }%
    }
    
    \DeclareFieldFormat{title}{\usebibmacro{string+doiurlisbn}{\mkbibemph{#1}}}
    \DeclareFieldFormat[article,incollection,thesis,inbook]{title}%
        {\usebibmacro{string+doiurlisbn}{\mkbibquote{#1}}}
    
	\usepackage{hyperref}
    \hypersetup{colorlinks=true,allcolors=blue}
    \usepackage[capitalize,nameinlink,noabbrev]{cleveref}
    \crefname{appsec}{Appendix}{Appendices}

    \usepackage{fancyhdr,lastpage}
    \theoremstyle{plain}
    \newtheorem{theorem}{Theorem}[section]
    \newtheorem{lemma}[theorem]{Lemma}
    
    \newtheorem{proposition}[theorem]{Proposition}
    
    \theoremstyle{definition}
    \newtheorem{definition}[theorem]{Definition}
    
    \theoremstyle{remark}
    
    \newtheorem{remark}[theorem]{Remark}
    
    \newtheorem{notation}[theorem]{Notation}
    \newtheorem{question}[theorem]{Question}

    \newcommand*\icon[1]{\tikz[baseline=(char.base)]{\node[shape=circle,draw,inner sep=1pt, color=black,text=black] (char) {#1};}}
    \usepackage{marvosym}

        \newcommand\restr[2]{{\left.\kern-\nulldelimiterspace #1 \vphantom{\big|} \right|_{#2} }}

    \title{Mahler's measure and elliptic curves \\ with potential complex multiplication}
    \author{Riccardo Pengo}
    \date{}

    \fancypagestyle{firstpage}
        {
            \fancyhf{}

            \lfoot{\small\icon{\Letter} \footnotesize Department of Mathematical Sciences, University of Copenhagen, Universitetsparken 5, 2100 Copenhagen, DENMARK \newline \small\icon{\Email} \footnotesize\href{mailto:pengo@math.ku.dk}{pengo@math.ku.dk}
            \small\icon{\Mundus} \footnotesize \href{https://sites.google.com/view/riccardopengo/}{https://sites.google.com/view/riccardopengo/}
            \small\icon{\Info} \footnotesize \href{https://zbmath.org/static/msc2020.pdf}{MSC 2020}: 11R06, 19F27, 14K22, 11S40}
        }

\begin{document}
\maketitle
\thispagestyle{firstpage}

\begin{abstract}
    Given an elliptic curve $E$ defined over $\mathbb{Q}$ which has potential complex multiplication by the ring of integers $\mathcal{O}_K$ of an imaginary quadratic field $K$ we construct a polynomial $P_E \in \mathbb{Z}[x,y]$ which is a planar model of $E$ and such that the Mahler measure $m(P_E) \in \mathbb{R}$ is related to the special value of the $L$-function $L(E,s)$ at $s = 2$.
\end{abstract}

\section{Introduction}

Let $E$ be an elliptic curve defined over a number field $F$ and let $L(E,s)$ be its $L$-function (see \cite[\S~C.16]{si09} for an introduction). 
We know that if $E$ has \textit{potential complex multiplication} (\textit{i.e.} $\operatorname{End}_{\overline{F}}(E) \not\cong \mathbb{Z}$) then $L(E,s)$ is an entire function, defined over the whole complex plane $\mathbb{C}$, which satisfies a functional equation relating $L(E,s)$ to $L(E,2-s)$ (in the ``arithmetic'' normalisation, that we will use in this paper). For example when $F = \mathbb{Q}$ this functional equation implies that $4 \pi^2 \, L'(E,0) = \mathfrak{f}_E \, L(E,2) \in \mathbb{R}_{> 0}$, where $\mathfrak{f}_E \in \mathbb{N}$ denotes the conductor of $E$ (see \cite[\S~IV.10]{si94}).

The aim of this paper is to prove that the \textit{special value} $L'(E,0)$ of the $L$-function associated to an elliptic curve $E$ defined over $\mathbb{Q}$ which has potential complex multiplication can be related to the \textit{Mahler measure} 
\[
    m(P) := \int_0^1 \int_0^1 \log\lvert P(e^{2 \pi i \theta_1},e^{2 \pi i \theta_2}) \rvert \, d\theta_1 \, d\theta_2 \in \mathbb{R}
\]
of a planar model $P \in \mathbb{Z}[x,y]$ of the elliptic curve $E$, as in the following theorem.

\begin{theorem}[see \cref{thm:main_theorem_body}] \label{thm:main_theorem}
    Let $E$ be an elliptic curve defined over $\mathbb{Q}$ such that $\operatorname{End}(E_{\overline{\mathbb{Q}}})\cong \mathcal{O}_K$ for some imaginary quadratic field $K$. Then there exists a polynomial $P \in \mathbb{Z}[x,y]$ such that: 
    \begin{itemize}
        \item its zero locus $V_P \hookrightarrow \mathbb{G}_m^2$ is birationally equivalent to $E$;
        \item $m(P) = r \, L'(E,0) + \log\lvert s \rvert$ for two explicit numbers $r \in \mathbb{Q}^{\times}$ and $s \in \overline{\mathbb{Q}}^{\times}$ defined in \eqref{eq:explicit_numbers}. 
    \end{itemize}
\end{theorem}

\cref{thm:main_theorem} fits into the vast landscape of \textit{Boyd's conjectures} (see \cref{que:boyd_conjectures}), which are precisely concerned with the relations between special values of $L$-functions and Mahler measures of polynomials (see \cref{sec:historical_introduction} for a brief historical review). 

Before moving on, two remarks about \cref{thm:main_theorem} are in order.
\begin{remark}
    The polynomial appearing in \cref{thm:main_theorem} will in general not be \textit{tempered} (see \cite[Section~2]{bo98} and \cite[Section~8]{vi99}). 
    Tempered polynomials have been traditionally the main focus of research on Boyd's conjectures, because for these polynomials the coordinate symbol $\{x,y\} \in H^2_{\mathcal{M}}(V_P,\mathbb{Q}(2))$ extends to an element of $H^2_{\mathcal{M}}(\widetilde{V_P},\mathbb{Q}(2))$, where $\widetilde{V_P}$ denotes a desingularisation of a compactification $\overline{V_P}$ of $V_P$.
    
    Nevertheless, in recent years more and more attention has been given to Mahler measures of non-tempered polynomials (see \cite{lsz16},\cite{lm18}, \cite{Meemark_Samart_2019}, \cite{gi19}) and our \cref{thm:main_theorem} fits into this history of examples.
    More precisely, in our case $\widetilde{V_P}$ is an elliptic curve and the zeros and poles of $x$ and $y$ are torsion points. 
    This helps to obtain a relation between $m(P)$ and $L'(E,0)$ despite the fact that $P$ is not tempered, because it allows to find a symbol $\eta_{x,y} \in H^2_{\mathcal{M}}(\widetilde{V_P},\mathbb{Q}(2))$ which is closely related to $\{x,y\} \in H^2_{\mathcal{M}}(V_P,\mathbb{Q}(2))$.
\end{remark}
\begin{remark}
    The polynomial appearing in \cref{thm:main_theorem} will have in general a very high degree, and thus the curve $V_P \subseteq \mathbb{G}_m^2$ will be highly singular. 
    This is in contrast with the majority of previously known cases of Boyd's conjectures, where the polynomials appearing have small degree.
\end{remark}

\subsection{Motivation and historical remarks}
\label{sec:historical_introduction}
$L$-functions and height theory have been two leading subjects in arithmetic geometry: the first ones are analytic counterparts of arithmetic and automorphic objects, and the second ones provide a way of measuring the complexity of these geometric objects.

Most $L$-functions are defined using an \textit{Euler product} (\textit{i.e.} a product indexed by prime ideals or closed points) which converges only for some complex numbers $s \in \mathbb{C}$, typically the ones such that $\Re(s) > \sigma$ for some $\sigma \geq 1$. It is expected that such an Euler product should be analytically continued to a meromorphic function belonging to the \textit{Selberg class} $\mathcal{S}$ (see \cite{Perelli_2005}).
This class consists of functions $f \in \mathcal{S}$ which admit an expansion as a Dirichlet series and as an Euler product in some half-plane, which are almost entire (\textit{i.e.} $(s-1)^m \, f(s)$ is entire for some $m \in \mathbb{N}$) and satisfy a suitable functional equation and a suitable ``growth condition'' (known as \textit{Ramanujan hypothesis}).

Functions belonging to $\mathcal{S}$ are known to enjoy a variety of properties, including the fact that they can be recovered from the collection of their \textit{special values} at the integers.
More precisely, if $f \colon \mathbb{C} \to \mathbb{C}$ is a meromorphic function we define its special value at $s_0 \in \mathbb{C}$ as
\begin{equation} \label{eq:special_values}
    f^{\ast}(s_0) := \lim_{s \to s_0} (s - s_0)^{-\operatorname{ord}_{s_0}(f)} \, f(s) \in \mathbb{C}^{\times}
\end{equation} 
where $\operatorname{ord}_{s_0}(f) \in \mathbb{Z}$ is the unique integer such that the limit \eqref{eq:special_values} exists and is different from zero.
We have then the following result, which is due to Deninger.

\begin{theorem}[see {\cite[Theorem~2.1]{de00}}] \label{thm:deninger_reconstruction}
    There exists a class of holomorphic functions $\mathcal{C}$ such that:
    \begin{itemize}
        \item for every $f \in \mathcal{S}$ and every $m \in \mathbb{N}$ such that $(s-1)^m \, f(s)$ is holomorphic we have that $(s-1)^m \, f(s) \in \mathcal{C}$;
        \item for every $n_0 \in \mathbb{Z}$ the map $\mathcal{C} \to (\mathbb{C}^{\times})^{\mathbb{Z}_{\geq n_0}}$ defined by $f \mapsto \{ f^{\ast}(n) \}_{n \geq n_0}$ is injective.
    \end{itemize}
\end{theorem}

\cref{thm:deninger_reconstruction} provides a very strong reason to study special values of functions in the Selberg class. Since the functions in the Selberg class satisfy a functional equation, we may restrict ourselves to the study of special values at \textit{non-positive integers}, for which we usually have cleaner formulas.

\begin{remark}
    Functions in the Selberg class are normalised ``analytically'', \textit{i.e.} they satisfy a functional equation which relates $s$ to $1 - s$. 
    Another possible way of normalising $L$-functions associated to arithmetic or automorphic object implies that $L(s)$ is related by the (conjectural) functional equation to $L^{\vee}(1-s)$, for some ``dual'' $L$-function $L^{\vee}$. 
    
    This often implies that $L(s)$ is related to $L(w - s)$ for some integer $w \in \mathbb{N}$, \textit{e.g.} $L(E,s)$ is related to $L(E,2-s)$ for the $L$-function associated to an elliptic curve.
    In the rest of this paper we will use this ``arithmetic'' normalisation of $L$-functions. This does not alter the conclusion of \cref{thm:deninger_reconstruction}, \textit{i.e.} $L$-functions which are normalised arithmetically can still be reconstructed from their special values.
\end{remark}

The study of special values of $L$-functions which come from arithmetic objects (that are conjectured to belong to the Selberg class) was initiated by Euler's proof of the identity $k \, \zeta(1-k) = - B_k$ for every $k \in \mathbb{N}_{\geq 1}$ (see \cite[Theorem~VII.1.8]{ne99}) and by Dirichlet's proof of the (analytic) class number formula (see \cite[Corollary~VII.5.11]{ne99}).
This formula relates the special value $\zeta_K^{\ast}(0)$ to the \textit{regulator} $R_K \in \mathbb{R}$, which can be seen as a height measuring the complexity of the number field $K$ (see \cite{Pazuki_2014}).

Other examples of links between values of $L$-functions and heights are given by the conjectural formulas of Birch and Swinnerton-Dyer (see \cite{Tate_1966}) and Colmez (see \cite{co93}).
These formulas fit into bigger landscapes of conjectures on special values of $L$-functions:
\begin{itemize}
    \item the formula of Birch and Swinnerton-Dyer is generalised by the Tamagawa number conjecture of Bloch and Kato (see \cite{bk90}, \cite{fpr94}), which followed work of Beilinson (see \cite{be84});
    \item the formula of Colmez has been recently generalised by Maillot and R{\"o}ssler (see \cite{Maillot_Rossler_2018}).
\end{itemize}

This paper studies another type of height which is conjectured to be related to special values of $L$-functions, which is defined as follows.
\begin{definition}
    The \textit{Mahler measure} of a Laurent polynomial
    $P \in \mathbb{C}[x_1^{\pm 1},\dots,x_n^{\pm 1}] \setminus \{ 0 \}$ is defined as
    \begin{equation} \label{eq:mahler_measure_definition}
        m(P) := \int_{\mathbb{T}^n} \log\lvert P \rvert \, d\mu_{\mathbb{T}^n} = \int_0^1 \cdots \int_0^1 \log\lvert P(e^{2 \pi i t_1},\dots,e^{2 \pi i t_n}) \rvert d t_1 \dots d t_n
    \end{equation}
    where $\mathbb{T}^n := (S^1)^n$ is the real $n$-torus and $\mu_{\mathbb{T}^n} := \frac{1}{(2 \pi i)^n} \, \frac{d z_1}{z_1} \wedge \dots \wedge \frac{d z_n}{z_n}$ is the unique Haar measure on $\mathbb{T}^n$ such that $\mu_{\mathbb{T}^n}(\mathbb{T}^n) = 1$.
\end{definition}

The integral appearing in \eqref{eq:mahler_measure_definition} is always convergent, $m(P) \in \mathbb{R}_{\geq 0}$ whenever $P \in \mathbb{Z}[x_1^{\pm 1},\dots,x_n^{\pm 1}]$ and there is an explicit classification of all the polynomials such that $m(P) = 0$ (see \cite[Chapter~3]{Everest_Ward_1999}).
Moreover, Mahler's measure of polynomials in one variable is linked to the logarithmic Weil height of algebraic numbers (see \cite[Proposition~1.6.6]{Bombieri_Gubler_2001}) and therefore to questions of Diophantine nature, such as \textit{Lehmer's problem}. 
This asks whether or not the set $\{ m(P) \colon P \in \mathbb{Z}[x] \} \setminus \{ 0 \} \subseteq \mathbb{R}_{> 0}$ has a minimum, and it seems that it could be approached by studying the Mahler measure of polynomials in multiple variables (see \cite{Smyth_2018}).
These Mahler measures appear to be far more mysterious than their one-variable counterparts: in particular, they appear to be related to special values of $L$-functions.

The two initial examples of these relationships are given by
\begin{equation} \label{eq:smyth_results}
    \begin{aligned}
        m(x+y+1) &= L^{\ast}(\chi_{-3},-1) \\
        m(x+y+z+1) &= -14 \, \zeta^{\ast}(-2)
    \end{aligned}
    \quad \text{where} \quad
    \begin{aligned}
        \chi_{-3} \colon (\mathbb{Z}/3 \mathbb{Z})^{\times} &\to \mathbb{C}^{\times} \\
        \pm 1 &\mapsto \pm 1
    \end{aligned}
\end{equation}
and have been proved by Smyth (see \cite{Smyth_1981}, \cite[Appendix~1]{bo81}). 
These initial successes led to an extensive search for more relations between Mahler measures and special values of $L$-functions, which culminated with the pioneering work \cite{bo98}. The numerical observations contained in Boyd's paper lead naturally to the following question.

\begin{question}[Boyd] \label{que:boyd_conjectures}
Let $P \in \mathbb{Z}[x_1^{\pm 1},\dots,x_n^{\pm 1}]$ be a Laurent polynomial such that $m(P) \neq 0$. 
Let $V_P \hookrightarrow \mathbb{G}_{m,\mathbb{Q}}^n$ be the zero locus of $P$, $\overline{V_P} \hookleftarrow V_P$ be a compactification of $V_P$ and $\widetilde{V_P} \twoheadrightarrow \overline{V_P}$ be a desingularisation of $\overline{V_P}$. 
Let $\ell \in \mathbb{Z}$ be a prime number and let $\operatorname{H}_{\ell}^{n - 1}(\widetilde{V_P}) := \operatorname{H}^{n - 1}_{\text{{\'e}t}}((\widetilde{V_P})_{\overline{\mathbb{Q}}};\mathbb{Q}_{\ell})$ denote the Galois representation given by the $\ell$-adic \'etale cohomology of $\widetilde{V_P}$ (see \cite[\href{https://stacks.math.columbia.edu/tag/03N1}{Chapter 03N1}]{stacks-project}). Let finally $r_P \in \mathbb{R}^{\times}$ be the real number
\[
r_P := \frac{L^{\ast}(\operatorname{H}_{\ell}^{n-1}(\widetilde{V_P}),0)}{m(P)}
\]
which is defined assuming that the Euler product defining $L(\operatorname{H}_{\ell}^{n-1}(\widetilde{V_P}),s)$ can be analytically continued until $s = 0$.
When is it true that $r_P \in \mathbb{Q}^{\times}$?
\end{question}

\begin{remark}
    Usual conjectures on $L$-functions imply that $L(\operatorname{H}^{n-1}_{\ell}(\widetilde{V_P}),s)$ should not depend on the choice of the prime $\ell$ (see \cite[Conjecture~1.3]{Taylor_2004}).
    Moreover, the special value $L^{\ast}(\operatorname{H}^{n-1}_{\ell}(\widetilde{V_P}),0)$ should be related to $L(\operatorname{H}^{n-1}_{\ell}(\widetilde{V_P}),n)$ by the functional equation (see \cite[\S~4.3.2]{fpr94}). 
    Hence \cref{que:boyd_conjectures} asks when the Mahler measure $m(P)$ of a polynomial $P \in \mathbb{Z}[x_1^{\pm 1},\dots,x_n^{\pm 1}]$ is related to special values of $L$-functions at $s = n$, where $n$ is the number of variables of $P$.
\end{remark}
\begin{remark}
\cref{que:boyd_conjectures} does not explain Smyth's results \eqref{eq:smyth_results}. 
These and other computations can be explained using an idea of Maillot (see \cite[\S~8]{Boyd_Lind_Rodriguez_Villegas_Deninger_2003}) concerning the \textit{reciprocal polynomial} $P^{\ast}(x_1,\dots,x_n) := P(x_1^{-1},\dots,x_n^{-1})$ (see \cite{la07} and the forthcoming work \cite{Brunault_Pengo}).
\end{remark}

Several works have answered positively to \cref{que:boyd_conjectures}, mostly using one of the following techniques:
\begin{itemize}
    \item functional equations for Mahler's measure of families of polynomials, which relate $m(P)$ to some kind of special function (often, hypergeometric) that is known to have a link to the special value (see \cite{Lalin_Rogers_2007}, \cite{Rogers_2011}, \cite{Rogers_Zudilin_2012}, \cite{Rogers_Zudilin_2014});
    \item explicit formulas for the regulator of modular functions (typically, modular units), which allow to relate it directly to the special value, and thus provide new ways of proving Beilinson's conjecture in specific cases (see \cite{Zudilin_2014}, \cite{Brunault_2016}, \cite{Brunault_2017});
    \item relations between different kinds of regulators, which allow to reduce some identities to other previously proved ones (see \cite{Lalin_Ramamonjisoa_2017}, \cite{Lalin_Wu_2018}).
\end{itemize}

For a survey of these techniques and for a history of Mahler's measure we highly recommend the forthcoming book \cite{Brunault_Zudilin_2020}.

\subsection{Our contribution and future research}
Let us explain how our main result \cref{thm:main_theorem} fits into the history outlined in \cref{sec:historical_introduction}. First of all, \cref{thm:main_theorem} provides a positive answer to \cref{que:boyd_conjectures} \textit{up to a logarithmic factor} for an infinite number of elliptic curves.
To be fair, we should observe that there are only finitely many elliptic curves with potential complex multiplication defined over $\mathbb{Q}$ \textit{up to twist}, \textit{i.e.} up to isomorphisms defined over $\overline{\mathbb{Q}}$. Hence \cref{thm:main_theorem} can be also seen as a step towards a positive answer to the following question.

\begin{question}[Twisting Boyd's conjectures]
    Let $P \in \mathbb{Z}[x_1^{\pm 1},\dots,x_n^{\pm 1}]$ be a Laurent polynomial such that $m(P) \neq 0$ and $L^{\ast}(\operatorname{H}_{\ell}^{n-1}(\widetilde{V_P}),0)/m(P) \in \mathbb{Q}^{\times}$. Let $X$ be a smooth, projective variety such that there exists an isomorphism $X_{\overline{\mathbb{Q}}} \xrightarrow{\sim} (\widetilde{V_P})_{\overline{\mathbb{Q}}}$. Does there exist a Laurent polynomial $Q \in \mathbb{Z}[x_1^{\pm 1},\dots,x_n^{\pm 1}]$ such that $V_Q$ is birational to $X$ and $L^{\ast}(\operatorname{H}_{\ell}^{n-1}(X),0)/m(Q) \in \mathbb{Q}^{\times}$?
\end{question}

Another element of novelty in this work is that we start from the elliptic curve $E$ and then subsequently look for a model of $E$ whose Mahler measure would be related to the special value $L'(E,0)$.
This differs from most of the past research, which usually starts from a polynomial $P$ (or from a family of polynomials $P_k$, such as $P_k \colon x + x^{-1} + y + y^{-1} + k$) and then uses the techniques outlined in \cref{sec:historical_introduction} to provide a link between the Mahler measure $m(P)$ (or the family of Mahler measures $m(P_k)$) and the special value $L^{\ast}(\operatorname{H}_{\ell}^{n-1}(\widetilde{V_P}),0)$ (or the family of special values $L^{\ast}(\operatorname{H}_{\ell}^{n-1}(\widetilde{V_{P_k}}),0)$).

Let us explain what is the strategy behind the proof of \cref{thm:main_theorem}. 
We know, thanks to the work of Deninger and Wingberg (see \cite{dw88}) and Rohrlich (see \cite{ro87}), that for every elliptic curve $E$ as in the statement of \cref{thm:main_theorem} there exist many pairs of functions $f,g \in \mathbb{Q}(E)$ such that the regulator of the \textit{Milnor symbol} $\{f,g\}$ is related to the special value $L'(E,0)$.
We prove in \cref{sec:defining_the_polynomial} that $\mathbb{Q}(E) = \mathbb{Q}(f,g)$, generalising a result of Brunault (see \cite{br16}).
This allows us to construct the polynomial $P \in \mathbb{Z}[x,y]$ as the minimal polynomial of $f$ and $g$. 
Finally, we can prove \cref{thm:main_theorem} by relating the regulator of $\{f, g\}$ to the Mahler measure of $P$, and this is done in \cref{sec:mahler_measure} using some generalisations of the seminal work of Deninger (see \cite{de97}).

Let us conclude this section with some questions which will serve as a guide for future research.

\begin{question} \label{que:removing_the_logarithm}
    Can we remove the logarithm from \cref{thm:main_theorem}? More precisely, given an elliptic curve $E$ defined over $\mathbb{Q}$ and such that $\operatorname{End}_{\overline{\mathbb{Q}}}(E) \cong \mathcal{O}_K$, does there exist a polynomial $P \in \mathbb{Z}[x,y]$ such that $V_P$ is birationally equivalent to $E$ and $m(P) = r \, L'(E,0)$ for some $r \in \mathbb{Q}^{\times}$?
\end{question}
\begin{remark}
    A positive answer to \cref{que:removing_the_logarithm} would be the first instance of a complete positive answer to \cref{que:boyd_conjectures} for an infinite family of $L$-functions.
    
    We believe that the great freedom allowed in the choice of the two functions $f,g \in \mathbb{Q}(E)$ which give rise to the polynomial $P$ in \cref{thm:main_theorem} will enable us to answer \cref{que:removing_the_logarithm} in the positive.
\end{remark}

\begin{question}
    Can we simplify the polynomials appearing in \cref{thm:main_theorem} to get something of smaller degree?
\end{question}

\begin{question} \label{que:higher_genus}
    Can the techniques of this paper be generalised to polynomials $P \in \mathbb{Z}[x,y]$ such that $\tilde{V_P}$ has higher genus?
\end{question}
\begin{remark}
    The most natural way to approach \cref{que:higher_genus} would be to prove a higher-dimensional analogue of Rohrlich's result (see \cref{thm:rohrlich_general_theorem}).
    This has been recently achieved for Siegel varieties (using modular techniques, and without any mention of complex multiplication) in \cite{Lemma_2017} and \cite{Cauchi_Lemma_Jacinto_2019}. 
\end{remark}

\subsection{Outline of the paper}
To conclude this introductory section, let us give an outline of the paper: \cref{sec:preliminaries} recalls some necessary preliminaries, \cref{sec:defining_the_polynomial} constructs the minimal polynomial $P \in \mathbb{Z}[x,y]$ associated to a CM elliptic curve $E$ and \cref{sec:mahler_measure} contains the computation of the Mahler measure $m(P)$. We also include \hyperref[sec:conductors]{an appendix} as a reference for a proof of Deuring concerning the relation between the conductor of a CM elliptic curve and the conductor of the associated Hecke character.

\section{Preliminaries} \label{sec:preliminaries}

The aim of this section is to provide various preliminaries that will be needed in the rest of the paper, concerning elliptic curves with complex multiplication (see \cref{sec:CM_elliptic_curves}), motivic cohomology (see \cref{sec:motivic_cohomology_and_tame_symbols}) and Deligne cohomology of curves (see \cref{sec:Deligne_cohomology}).

\begin{notation}
    For us $\mathbb{N}$ denotes the monoid of natural numbers, so $0 \in \mathbb{N}$. Moreover, for every $a \in \mathbb{Z}$ we write $\mathbb{Z}_{\geq a} := \{n \in \mathbb{Z} \, \colon \, n \geq a \}$.
\end{notation}

\subsection{Elliptic curves with complex multiplication} \label{sec:CM_elliptic_curves}

We will use this section to recall some useful facts about elliptic curves with complex multiplication, following mainly \cite[Chapter~II]{si94} and \cite[Chapter~10]{la87}. We will also recall some facts and definitions coming from \cite{ro87}.

Let $E$ be an elliptic curve defined over $\mathbb{Q}$ which has \textit{potential complex multiplication} in the sense of \cref{def:complex_multiplication}. Then there exists a unique imaginary quadratic field $K$, a unique order $\mathcal{O} \subseteq \mathcal{O}_K$ and a unique isomorphism 
\[
    [\,]_E \colon \mathcal{O} \xrightarrow{\sim} \operatorname{End}(E_{\overline{\mathbb{Q}}})     
\]
such that $[\alpha]_E^{\ast}(\omega_E) = \alpha \, \omega_E$ for every $\alpha \in \mathcal{O}$, where $\omega_E \in \Omega^1_{E}$ denotes the invariant differential of $E$ (see \cite[Proposition~1.1]{si94}).
Moreover, we know that $\mathcal{O}$ belongs to the finite list of imaginary quadratic orders such that $\operatorname{Pic}(\mathcal{O}) = 1$ (see  \cite[Theorem~7.30]{co13} and \cite[Theorem~11.1]{co13}).

\begin{remark}
    In this paper we will only consider elliptic curves $E$ defined over $\mathbb{Q}$ which have potential complex multiplication by \textit{the maximal order} $\mathcal{O}_K$ of an imaginary quadratic field $K$, \textit{i.e.} such that $\operatorname{End}(E_{\overline{\mathbb{Q}}}) \cong \mathcal{O}_K$. 
    We believe nevertheless that most of our results will hold in the general case. 
\end{remark}

\begin{notation} \label{not:CM_notation_1}
    Let $E$ be an elliptic curve defined over $\mathbb{Q}$, having potential complex multiplication by the ring of integers $\mathcal{O}_K$ of an imaginary quadratic field $K$. We introduce the following notation: 
    \begin{itemize}
        \item $\mu_K := \mathcal{O}_K^{\times}$ is the group of units, which is finite;
        \item $E(\mathbb{R})^0 \subseteq E(\mathbb{C})$ denotes the connected component of the identity of the real Lie group $E(\mathbb{R})$;
        \item $\gamma(E_{\mathbb{R}}/\mathbb{R}) \in \{\pm 1\}$ is defined as
            \[
                \gamma(E_{\mathbb{R}}/\mathbb{R}) := \begin{cases} - \operatorname{sign}(b), \ \text{if} \ b \neq 0 \ (\text{\textit{i.e.}} j(E) \neq 1728) \\ - \operatorname{sign}(a), \ \text{if} \ b = 0 \ (\text{\textit{i.e.}} j(E) = 1728)  \end{cases}
            \]
        for any short Weierstrass equation $y^2 = x^3 + a \, x + b$ defining $E$ (see \cite[Proposition~V.2.2]{si94});
        \item $\mathfrak{b}_E \subseteq K$ is a fractional ideal defined as
            \[
                \mathfrak{b}_E :=
                \begin{cases} \mathcal{O}_K, \ \text{if} \ (\operatorname{disc}(K/\mathbb{Q}) \neq -4 \land \gamma(E_{\mathbb{R}}/\mathbb{R}) = -1) \lor (\operatorname{disc}(K/\mathbb{Q}) = -4 \land 4 a^3 + 27 b^2 < 0) \\ 
                \mathfrak{D}_{K/\mathbb{Q}}^{-1}, \ \text{otherwise}  \end{cases}
            \]
        for any short Weierstrass equation $y^2 = x^3 + a \, x + b$ defining $E$ (recall that the discriminant of such a Weierstrass equation is defined as $-16 (4 a^3 + 27 b^2)$). 
        Moreover, 
        \[
            \mathfrak{D}_{K/\mathbb{Q}}^{-1} := \{ x \in K \, \mid \, \operatorname{Tr}(x \, \mathcal{O}_K) \subseteq \mathbb{Z} \}
        \]
        denotes the inverse of the different ideal associated to the extension $\mathbb{Q} \subseteq K$ (see \cite[Section~III.2]{ne99});
        \item $b_E \in K$ is any generator of $\mathfrak{b}_E$ as an $\mathcal{O}_K$-module;
    \end{itemize}
\end{notation}

We now recall the existence of a particular complex uniformisation $\theta_E \colon \mathbb{C} \twoheadrightarrow E(\mathbb{C})$ for a CM elliptic curve $E$ defined over $\mathbb{Q}$, which is due to Rohrlich (see \cite[Page~377]{ro87}).

\begin{proposition}[Rohrlich] \label{prop:rohrlich_uniformization}
    For every embedding $\sigma \colon K \hookrightarrow \mathbb{C}$ and every orientation of $E(\mathbb{R})^{0}$ there exists a unique surjective map of complex Lie groups $\theta_E \colon \mathbb{C} \twoheadrightarrow E(\mathbb{C})$ such that $\ker(\theta_E) = \sigma(\mathfrak{b}_E)$, $\theta_E(\mathbb{R}) \subseteq E(\mathbb{R})$ and the induced isomorphism of real Lie groups $\mathbb{R}/\mathbb{Z} \xrightarrow{\sim} E(\mathbb{R})^0$ preserves the orientations.
\end{proposition}

\begin{remark}
    To avoid unnecessary sign issues, whenever we have an elliptic curve $E$ defined over $\mathbb{Q}$ which has potential complex multiplication we will fix implicitly an embedding $\sigma \colon K \hookrightarrow \mathbb{C}$ and an orientation of $E(\mathbb{R})^{0}$.
\end{remark}

Let us recall some more notation associated to an elliptic curve $E$ with potential complex multiplication. 

\begin{notation}
    Let $E$ be an elliptic curve defined over $\mathbb{Q}$, having potential complex multiplication by $\mathcal{O}_K$. We introduce the following notation:
    \begin{itemize}
        \item $\mathbb{A}_K$ denotes the ring of ad{\`e}les of $K$ (see \cite[Section~VI.1]{ne99});
        \item $\psi_E \colon \mathbb{A}_K^{\times}/K^{\times} \to K^{\times}$ denotes the algebraic Hecke character associated to $E$ (see \cite[Section~1.1]{sc88}).
        This Hecke character is related to $E$ by the fact that $L(E,s) = L(\sigma \circ \psi_{E},s)$ for any embedding $\sigma \colon K \hookrightarrow \mathbb{C}$ (see \cite[Page~187]{mi72});
        \item $\mathfrak{f}_{\psi_E} \subseteq \mathcal{O}_K$ denotes the conductor of $\psi_E$ (see \cite[Section~VII.6]{ne99});
        \item abusing notation, we denote also by $\psi_E \colon I(\mathfrak{f}_{\psi_E}) \to K^{\times}$ the classical algebraic Hecke character associated to $\psi_E$ (see \cite[Section~0.1]{sc88}). Here $I(\mathfrak{f}_{\psi_E})$ denotes the group of fractional ideals of $K$ which are coprime to $\mathfrak{f}_{\psi_E}$;
        \item $\chi_E \colon (\mathcal{O}_K/\mathfrak{f}_{\psi_E})^{\times} \to \mu_K$ is the unique group homomorphism such that $\psi_E(x \, \mathcal{O}_K) = \chi_E(x) \, \sigma(x)$ for every $x \in \mathcal{O}_K$ coprime with the conductor $\mathfrak{f}_{\psi_E}$. 
        Here, by a slight abuse of notation, we denote by $\sigma \in \operatorname{Gal}(K/\mathbb{Q})$ the element corresponding to the embedding $\sigma \colon K \hookrightarrow \mathbb{C}$.
        We extend $\chi_E$ by zero to get a multiplicative map $\chi_E \colon \mathcal{O}_K \to \mu_K \cup \{ 0 \}$, and we observe that $\overline{\chi_E(x)} = \chi_E(\overline{x})$ for every $x \in \mathcal{O}_K$;
        \item for every $\alpha \in \mathcal{O}_K$ such that $\chi_E(\alpha) = 0$ we can define a map 
        \begin{align*}
            \widetilde{\chi_E} \colon E[\alpha](\mathbb{C}) &\to \mu_K \cup \{ 0 \} \\
            x &\mapsto \chi_E\left( \theta_E^{-1}(x) \, \frac{\overline{\alpha}}{b_E} \right)
        \end{align*}
        where $E[\alpha] := \ker(E \xrightarrow{[\alpha]_E} E)$ denotes the group (scheme) of $\alpha$-torsion points. 
        We observe that if $\overline{\alpha}/b_E \in \mathbb{R}$ (\textit{e.g.} if $\alpha = \nu_E$) then $\overline{\widetilde{\chi_E}(x)} = \widetilde{\chi_E}(\overline{x})$ for every $x \in E[\alpha](\mathbb{C})$;
        \item by abuse of notation, for every $\alpha \in \mathcal{O}$ we write $E[\alpha]$ for the full group $E[\alpha](\mathbb{C})$ of $\alpha$-torsion points. Moreover, we denote by $K(E[\alpha])$ the smallest extension of $K$ over which all the $\alpha$-torsion points are defined;
        \item $\nu_E \in \mathcal{O}_K$ is defined as $\nu_E := \operatorname{N}_{K/\mathbb{Q}}(\mathfrak{f}_{\psi_{E}}) \, \min\{ \mathfrak{b}_E \cap \mathbb{R}_{> 0} \}/b_E$.
    \end{itemize}
\end{notation}

We will later use the group of $\nu_E$-torsion points $E[\nu_E]$ (see \cref{not:CM_notation_1}) to define a model for $E$. 
In particular we will need the following description of the Galois action on these points.

\begin{lemma}[see {\cite[Page~264]{dw88}}] \label{lem:galois_action_on_torsion}
    We can describe the action of $\operatorname{Gal}(K(E[\nu_E])/K)$ on $E[\nu_E]/\mu_K$ as
    \begin{align*}
        \operatorname{Gal}(K(E[\nu_E])/K) \times E[\nu_E]/\mu_K &\to E[\nu_E]/\mu_K \\
        (\sigma, x) &\mapsto [\varphi(\restr{\sigma}{R_{\nu_E}})^{-1} \, \widetilde{\chi_E}(x)]_E(x)
    \end{align*}
    where $R_{\nu_E} \subseteq K(E[\nu_E])$ denotes the ray class field of $K$ relative to $\nu_E \, \mathcal{O}_K$ (see \cite[Definition~VI.6.2]{ne99}) and 
    \[
        \varphi \colon \operatorname{Gal}(R_{\nu_E}/K) \xrightarrow{\sim} (\mathcal{O}_K/\nu_E)^{\times}/\mu_K
    \]
    denotes the isomorphism which follows from \cite[Theorem~II.5.6]{si94}.
\end{lemma}
\begin{remark}
    Since $E$ is clearly a module over $\operatorname{End}(E)$ we have that $E[\nu_E]$ is a module over $\mathcal{O}_K/\nu_E$, and thus it makes sense to act over $E[\nu_E]/\mu_K$ by $\varphi(\restr{\sigma}{R_{\nu_E}})^{-1} \in (\mathcal{O}_K/\nu_E)^{\times}/\mu_K$. Moreover, the action of $\widetilde{\chi_E}(x)$ over $E[\nu_E]/\mu_K$ is either given by the identity or by the zero map.
\end{remark}

Let us introduce some more notation, including the definition of the ``diamond'' operator $\diamondsuit$ and of the function $\mathcal{R}$ which appear in \cref{thm:main_theorem_body}. 

\begin{notation} \label{not:orders_and_divisors}
    Let $E$ be an elliptic curve over a field $\kappa$. We introduce the following notation:
    \begin{itemize}
        \item for every function $f \colon E \to \mathbb{P}^1$ we denote by $S_f \subseteq E(\overline{\kappa})$ the set of zeros and poles of $f$, and by $n_f \in \mathbb{Z}_{\geq 1} \cup \{\infty\}$ the least common multiple of the orders of the points in $S_f$;
        \item for every set of functions $f_1,\dots,f_k \colon E \to \mathbb{P}^1$ we define $S_{f_1,\dots,f_k} \subseteq E(\overline{\kappa})$ as $S_{f_1,\dots,f_k} := \bigcup_{i = 1}^k S_{f_i}$, and $n_{f_1,\dots,f_k} \in \mathbb{Z}_{\geq 1} \cup \{\infty\}$ as $n_{f_1,\dots,f_k} := \operatorname{lcm}(n_{f_1},\dots,n_{f_k})$;
        \item for every set of functions $f_1,\dots,f_k \colon E \to \mathbb{P}^1$ such that $S_{f_1,\dots,f_k} \subseteq E(\overline{\kappa})_{\text{tors}}$ (which is equivalent to say that $n_{f_1,\dots,f_k} \neq \infty$) and every point $x \in S_{f_1,\dots,f_k}$ we denote by $\varphi_{f_1,\dots,f_n}^{(x)} \colon E \to \mathbb{P}^1$ any function such that $\operatorname{div}(\varphi_{f_1,\dots,f_n}^{(x)}) = n_{f_1,\dots,f_k} \, ((x) - (0))$. This function exists, and is uniquely determined up to constants;
        \item for every divisor $D \in \mathbb{Q}[E(\overline{\kappa})_{\text{tors}}]$ we define $\operatorname{ord}(D) \in \mathbb{Z}_{\geq 1}$ as the smallest natural number $n \in \mathbb{Z}_{\geq 1}$ such that $n \, D$ is a principal divisor. 
        More concretely, we can write $D$ uniquely as $D = \frac{1}{a} \, \sum_{j = 1}^J b_j (P_j)$ for some $a \in \mathbb{Z}_{\geq 1}$ and $\{b_j\}_{j = 1}^J \subseteq \mathbb{Z}$ such that $\operatorname{gcd}(a,b_j) = 1$ for every $j \in \{1,\dots,J\}$.
        Then $\operatorname{ord}(D)$ equals the order of the point $\sum_{j = 1}^J b_j P_j \in E(\overline{\kappa})_{\text{tors}}$.
    \end{itemize}
\end{notation}

\begin{notation} \label{not:diamond_operator}
    We introduce the following notation:
    \begin{itemize}
        \item for every group $G$ acting on a set $S$ (from the left) we denote by $G \, s \subseteq S$ the orbit of any element $s \in S$, by $G \backslash S$ the set of orbits and by $\operatorname{Stab}_G(s) \leq G$ the stabiliser of any element $s \in S$;
        \item for every field $\kappa$ we denote by $G_{\kappa} := \operatorname{Gal}(\overline{\kappa}/\kappa)$ its absolute Galois group;
        \item we recall that for every $\kappa$-variety $V$ the set $V(\overline{\kappa})$ is endowed with an action of $G_{\kappa}$, and we denote by $[x]_{\kappa} := G_{\kappa} \, x$ the orbits of this action. Moreover, we endow $V(\overline{\kappa})^2 := V(\overline{\kappa}) \times V(\overline{\kappa})$ with the diagonal action of $G_{\kappa}$, and we denote by $[x,y]_{\kappa}$ the orbits under this diagonal action;
        \item for every (smooth, projective) curve defined over a field $\kappa$ we denote by $J_X$ the Jacobian of $X$, with a chosen inclusion $\iota \colon X \hookrightarrow J_X$ (see \cite{Milne_86}). Recall that for every $x,y \in X(\overline{\kappa})$ the difference $x - y := \iota(x) - \iota(y) \in J_X(\overline{\kappa})$ does not depend on $\iota$;
        \item we observe that for every function $f \in \kappa(X)$ its divisor $\operatorname{div}(f) \in \mathbb{Z}[ X(\overline{\kappa})]$ is invariant under the action of $\mathcal{G}_{\kappa}$. This is to say that if $\operatorname{div}(f) = \sum_{x \in X(\overline{k})} a_x \, (x)$ then $a_{x} = a_{\sigma(x)}$ for every $\sigma \in G_{\kappa}$, hence we can write $a_{[x]_{\kappa}} := a_x$ for every Galois orbit $[x]_{\kappa} \in G_{\kappa} \backslash X(\overline{\kappa})$;
        \item finally, we define the diamond operator $\diamondsuit \colon (\mathbb{Z}[X(\overline{\kappa})]^{G_{\kappa}})^{\otimes 2} \to \mathbb{Z}[G_{\kappa} \backslash J_X(\overline{\kappa})]$ by
        \[
            \left( \sum_{x \in X(\overline{\kappa})} a_{x} \, (x) \right) \diamondsuit \left( \sum_{y \in X(\overline{\kappa})} b_{y} \, (y) \right) := \sum_{[x,y]_{\kappa} \in G_{\kappa} \backslash X(\overline{\kappa})^2} a_{[x]_{\kappa}} \, b_{[y]_{\kappa}} \, \frac{\lvert [x,y]_{\kappa} \rvert }{\lvert [x - y]_{\kappa} \rvert} \, ( [x - y]_{\kappa} )
        \]
        which coincides with the diamond operator defined in the previous literature (see \cite{Lalin_Ramamonjisoa_2017} and \cite{lm18}) for divisors supported on rational points. 
    \end{itemize}
\end{notation}

\begin{notation} \label{not:rational_factor_function}
    Let $E$ be an elliptic curve defined over $\mathbb{Q}$, having potential complex multiplication by $\mathcal{O}_K$. We introduce the following notation:
    \begin{itemize}
        \item for every point $x \in E(\overline{\mathbb{Q}})$ we define $\operatorname{Ann}_{\mathcal{O}_K}(x) := \{ \alpha \in \mathcal{O}_K \mid [\alpha](x) = 0 \}$, which is an $\mathcal{O}_K$-ideal;
        \item $\xi_x \in \mathcal{O}_K$ is the unique generator of the ideal $\operatorname{Ann}_{\mathcal{O}_K}(x) \, \mathfrak{b}_E^{-1} \subseteq \mathcal{O}_K$ such that $\theta_E(\xi_x^{-1}) \in G_K \, x$;
        \item we define a function
        \begin{align*}
            \mathfrak{s} \colon E(\overline{\mathbb{Q}})_{\text{tors}} &\to K \\
            x &\mapsto \frac{1}{\xi_x} \, \prod_{\mathfrak{p} \mid \operatorname{Ann}_{\mathcal{O}_K}(x)} \left( 1 - \psi_E(\mathfrak{p}) \right)
        \end{align*}
        where $\mathfrak{p}$ runs over all the prime ideals of $\mathcal{O}_K$ which divide $\operatorname{Ann}_{\mathcal{O}_K}(x)$;
        \item we define another function
        \begin{align*}
            \mathcal{R} \colon E(\overline{\mathbb{Q}})_{\text{tors}} &\to \mathbb{Q} \\
            x &\mapsto \begin{cases} 0, \ \text{if} \ \mathfrak{f}_{\psi_E} \nmid \operatorname{Ann}_{\mathcal{O}_K}(x) \\
            \mathfrak{s}(x), \ \text{if} \ \mathfrak{f}_{\psi_E} \mid \operatorname{Ann}_{\mathcal{O}_K}(x) \ \text{and} \ G_K \, x = G_K \, \overline{x} \\
            \mathfrak{s}(x) + \overline{\mathfrak{s}(x)}, \ \text{if} \ \mathfrak{f}_{\psi_E} \mid \operatorname{Ann}_{\mathcal{O}_K}(x) \ \text{and} \ G_K \, x \neq G_K \, \overline{x} \\
            \end{cases}
        \end{align*}
        and we observe that for every $x \in E(\overline{\mathbb{Q}})_{\text{tors}}$ we have that $\operatorname{ord}(x) \, \mathcal{R}(x) \in \mathbb{Z}$, that $\mathcal{R}([-1](x)) = - \mathcal{R}(x)$ and that $\mathcal{R}(\sigma(x)) = \mathcal{R}(x)$ for every $\sigma \in G_{\kappa}$; 
        \item we write $\mathcal{R} \colon \mathbb{Q}[E(\overline{\mathbb{Q}})_{\text{tors}}] \to \mathbb{Q}$ for the $\mathbb{Q}$-linear extension of $\mathcal{R} \colon E(\overline{\mathbb{Q}})_{\text{tors}} \to \mathbb{Q}$, and we observe that it descends to a $\mathbb{Q}$-linear map $\mathcal{R} \colon \mathbb{Q}[G_{\mathbb{Q}} \backslash E(\overline{\mathbb{Q}})_{\text{tors}}] \to \mathbb{Q}$.
    \end{itemize}
\end{notation}

Finally, let us recall a result proved recently by Campagna and the author which concerns \textit{division fields} of elliptic curves with complex multiplication. This will be used in \cref{sec:rohrlich_polynomials} to construct suitable functions on an elliptic curve with potential complex multiplication. 

\begin{proposition}[see {\cite[Corollary~1.3]{Campagna_Pengo}}]
\label{prop:campagna_pengo}    
    Let $E$ be an elliptic curve defined over $\mathbb{Q}$ having potential complex multiplication by the ring of integers $\mathcal{O}_K$ of some imaginary quadratic field $K$. 
    Let moreover $p \in \mathbb{N}$ be a prime number such that $p \, \mathcal{O}_K$ is a prime ideal (\emph{i.e.} $p$ is inert in $K$) and such that $p \nmid \mathfrak{f}_E$. Let finally $n \in \mathbb{N}$ be any natural number such that $p \nmid n$. 
    Then we have that $K(E[p]) \cap K(E[n]) = K$.
\end{proposition}

\subsection{Motivic cohomology and tame symbols} \label{sec:motivic_cohomology_and_tame_symbols}
Let $\kappa$ be a number field and $X$ be a curve over $\kappa$. The aim of this section is to recall the isomorphism
\begin{equation} \label{eq:tame_symbol}
    H^2_{\mathcal{M}}(X;\mathbb{Q}(2)) \cong \ker\left( 
    \begin{aligned}
        \frac{(\kappa(X)^{\times} \otimes_{\mathbb{Z}} \mathbb{Q})^{\otimes 2}}{\langle x \otimes (1 - x) \colon x \in \kappa(X)^{\times} \setminus \{ 1 \} \rangle} &\xrightarrow{\partial} \bigoplus_{x \in \lvert X \rvert} \kappa(x)^{\times} \otimes_{\mathbb{Z}} \mathbb{Q} \\
        \{f,g\} &\mapsto \bigoplus_{x \in \lvert X \rvert} \partial_x(\{f,g\})
    \end{aligned} \right)
\end{equation}
between the motivic cohomology $\mathbb{Q}$-vector space $H^2_{\mathcal{M}}(X;\mathbb{Q}(2))$ and the kernel of the \textit{tame symbol} map $\partial$, obtained by gluing the maps
\[
    \partial_x(\{f,g\}) := (-1)^{\operatorname{ord}_x(f) \operatorname{ord}_x(g)} \left. \frac{f^{\operatorname{ord}_x(g)}}{g^{\operatorname{ord}_x(f)}} \right\rvert_{x}
\]
where $\{f,g\}$ denotes the class of $f \otimes g$ in the quotient of $(\kappa(X)^{\times} \otimes_{\mathbb{Z}} \mathbb{Q})^{\otimes 2}$ by the ideal generated by the elements of the form $x \otimes (1 - x)$, with $x \in \kappa(X)^{\times} \setminus \{ 1 \}$.

First of all, recall that the bi-graded motivic cohomology groups $H^n_{\mathcal{M}}(S;\Lambda(m))$ (which appear in \eqref{eq:tame_symbol} for $n = m = 2$) can be defined as
\[
    H^n_{\mathcal{M}}(S;\Lambda(m)) := \operatorname{Hom}_{\mathrm{DM}(S,\Lambda)}(\mathbf{1}_S,\mathbf{1}_S(m)[n])
\]
for any ring $\Lambda$ and any scheme $S$ (see \cite[Section~11.2]{Cisinski_Deglise_2019}). Here, $\mathrm{DM}(S,\Lambda)$ denotes the triangulated category of mixed motives over $S$, as constructed in \cite{Cisinski_Deglise_2019}. This is defined as the $\mathbb{A}^1$-localisation of the derived category $\operatorname{D}(\mathrm{SH}^{\text{tr}}(S;\Lambda))$ of $\Lambda$-sheaves on $S$ which have transfers. It can be shown that there is a natural choice of tensor product which makes $\mathrm{DM}(S,\Lambda)$ into a monoidal category, and we denote by $\mathbf{1}_S \in \mathrm{DM}(S,\Lambda)$ the unit for this monoidal structure. Finally, $[n] \colon \mathrm{DM}(S,\Lambda) \to \mathrm{DM}(S,\Lambda)$ denotes the shift functor (coming from the fact that $\mathrm{DM}(S,\Lambda)$ is triangulated), and $(m) \colon \mathrm{DM}(S,\Lambda) \to \mathrm{DM}(S,\Lambda)$ denotes the Tate twist (see \cite[Section~1.1.d]{Cisinski_Deglise_2019}).

\begin{notation}
    From now on, let us use the notation $H^{n,m}_{\mathcal{M}}(S) := H^n_{\mathcal{M}}(S;\mathbb{Q}(m)) = \operatorname{Hom}_{\mathrm{DM}(S,\mathbb{Q})}(\mathbf{1}_S,\mathbf{1}_S(m)[n])$.
\end{notation}

\begin{remark} \label{rmk:properties_of_motivic_cohomology}
    Let us recall some properties of motivic cohomology, which will be useful in what follows:
    \begin{itemize}
        \item motivic cohomology commutes with disjoint unions. In particular if $Y$ is a discrete scheme defined over a field $\kappa$ then
        \begin{equation} \label{eq:motivic_cohomology_finitely_many_points}
            H^{n,m}_{\mathcal{M}}(Y) \cong \bigoplus_{y \in \lvert Y \rvert} H^{n,m}_{\mathcal{M}}(\operatorname{Spec}(\kappa(y)))
        \end{equation}
        where $\kappa(y) \supseteq \kappa$ denotes the residue field of $Y$ at $y$;
        \item for every regular scheme $S$ we have that $H^{1,1}_{\mathcal{M}}(S) \cong \mathcal{O}(S)^{\times} \otimes_{\mathbb{Z}} \mathbb{Q}$ (see \cite[Theorem~11.2.14]{Cisinski_Deglise_2019});
        \item we have that $H^{2 n,m}_{\mathcal{M}}(\operatorname{Spec}(F)) = 0$ for every number field $F$ and every $n, m \in \mathbb{Z}$. This follows from the identification between rational motivic cohomology and $K$-theory (see \cite[Corollary~14.2.14]{Cisinski_Deglise_2019}) and from Borel's theorem on the $K$-theory of number fields (see \cite{Soule_2010} for a survey).
    \end{itemize}
\end{remark}

Recall now that, if $\iota \colon W \hookrightarrow V$ is a closed immersion of pure co-dimension $c \in \mathbb{N}$ (see \cite[\href{https://stacks.math.columbia.edu/tag/04MS}{Section 04MS}]{stacks-project}) between schemes which are smooth, separated and of finite type over a common base $S$ then there is a long exact sequence (sometimes called \textit{localisation sequence} or \textit{Gysin sequence}) given by
\begin{equation} \label{eq:localisation_sequence}
    \begin{tikzcd}
         & \dots \rar \ar[draw=none]{d}[name=X, anchor=center]{} & H^{n - 1,m}_{\mathcal{M}}(V \setminus \iota(W))
        \ar[rounded corners, to path={ -- ([xshift=2ex]\tikztostart.east) |- (X.center) \tikztonodes -| ([xshift=-2ex]\tikztotarget.west) -- (\tikztotarget)}]{dll}[at end]{\delta} \\ %
        H^{n - 2 c,m-c}_{\mathcal{M}}(W) \arrow[r,"\iota_{\ast}"] & H^{n,m}_{\mathcal{M}}(V) \arrow[r,"j^{\ast}"] \ar[draw=none]{d}[name=X, anchor=center]{} & H^{n,m}_{\mathcal{M}}(V \setminus \iota(W)) \ar[rounded corners, to path={ -- ([xshift=2ex]\tikztostart.east) |- (X.center) \tikztonodes -| ([xshift=-2ex]\tikztotarget.west) -- (\tikztotarget)}]{dll}[at end]{\delta} \\
        H^{n + 1 - 2 c,m-c}_{\mathcal{M}}(W) \rar & \dots &
    \end{tikzcd}
\end{equation}
where $j \colon V \setminus \iota(W) \hookrightarrow V$ denotes the open embedding which is complementary to $\iota$ (see \cite[Section~11.3.4]{Cisinski_Deglise_2019}). 

Now, let us take $V = X$ to be a regular, connected curve (defined over a number field $\kappa$) and $n = m = 2$. 
Then if we use the properties of motivic cohomology recalled in \cref{rmk:properties_of_motivic_cohomology} in the localisation sequence \eqref{eq:localisation_sequence} we get an exact sequence
\begin{equation} \label{eq:short_localisation_sequence}
    0 \to H^{2,2}_{\mathcal{M}}(X) \to H^{2,2}_{\mathcal{M}}(X \setminus Y) \xrightarrow{\delta} H^{1,1}_{\mathcal{M}}(Y) \to H^{3,2}_{\mathcal{M}}(X) \to H^{3,2}_{\mathcal{M}}(X \setminus Y) \to 0
\end{equation}
for every finite set of closed points $Y \subseteq X$.
Moreover, we know again from \cref{rmk:properties_of_motivic_cohomology} that
\begin{equation} \label{eq:motivic_cohomology_divisors}
    H^{1,1}_{\mathcal{M}}(Y) \overset{\text{\eqref{eq:motivic_cohomology_finitely_many_points}}}{\cong} \bigoplus_{y \in \lvert Y \rvert} H^{1,1}_{\mathcal{M}}(\operatorname{Spec}(\kappa(y))) \cong \bigoplus_{y \in \lvert Y \rvert} \kappa(y)^{\times} \otimes_{\mathbb{Z}} \mathbb{Q}
\end{equation}
which shows that if we let $Y$ grow we get an exact sequence
\[
    0 \to H^{2,2}_{\mathcal{M}}(X) \to H^{2,2}_{\mathcal{M}}(\xi_X) \xrightarrow{\delta} \bigoplus_{x \in \lvert X \rvert} \kappa(x)^{\times} \otimes_{\mathbb{Z}} \mathbb{Q} \to \dots
\]
where $\xi_X \in X$ denotes the generic point. 
We can now get the isomorphism \eqref{eq:tame_symbol} using the identification
\[
    H^{2,2}_{\mathcal{M}}(\xi_X) \cong H^{2,2}_{\mathcal{M}}(\operatorname{Spec}(\kappa(X))) \cong \frac{(\kappa(X)^{\times} \otimes_{\mathbb{Z}} \mathbb{Q})^{\otimes 2}}{\langle x \otimes (1 - x) \colon x \not\in \{0,1\} \rangle}
\]
and its compatibility with localisation sequences in motivic cohomology and $K$-theory, which identifies $\partial$ with $\delta$.

Now let us use what we have just recalled to construct elements in motivic cohomology.
\begin{proposition} \label{prop:motivic_cohomology_decomposition}
    Let $X$ be a curve over a number field $\kappa$ and let $Y \hookrightarrow X$ be a finite set of closed points. Assume that there exists $y_0 \in Y$ such that $[\kappa(y_0) \colon \kappa] \, y - [\kappa(y) \colon \kappa] y_0 \in (J_X)_{\text{tors}}$ for every $y \in Y$. Then the natural restriction map $H^{2,2}_{\mathcal{M}}(X) \to H^{2,2}_{\mathcal{M}}(X \setminus Y)$ admits a natural retraction $H^{2,2}_{\mathcal{M}}(X \setminus Y) \twoheadrightarrow H^{2,2}_{\mathcal{M}}(X)$.
\end{proposition}
\begin{proof}
    Let $\kappa' \supseteq \kappa$ be a finite Galois extension, such that all the points of $Y$ are $\kappa'$-rational.
    Then the identification \eqref{eq:motivic_cohomology_divisors} gives the isomorphism 
        \begin{equation} \label{eq:motivic_cohomology_finite_set_base_change}
            H^{1,1}_{\mathcal{M}}(Y_{\kappa'}) \cong \bigoplus_{y \in \lvert Y \rvert} (\kappa')^{\times} \otimes_{\mathbb{Z}} \mathbb{Q} \cong \mathbb{Q}[Y] \otimes_{\mathbb{Z}} (\kappa')^{\times}    
        \end{equation}
    where $\mathbb{Q}[Y]$ denotes the group of divisors with rational coefficients which are supported on $Y$. 
    Now the exact sequence \eqref{eq:short_localisation_sequence} induces a short exact sequence
    \[
        0 \to H^{2,2}_{\mathcal{M}}(X_{\kappa'}) \to H^{2,2}_{\mathcal{M}}((X \setminus Y)_{\kappa'}) \xrightarrow{\delta'} \operatorname{Im}(\delta') \to 0
    \]
    and using Weil's reciprocity law (see \cite[\S 6.12.1]{Weibel_2013}) we can see that, under the isomorphism \eqref{eq:motivic_cohomology_finite_set_base_change} we have that $\operatorname{Im}(\delta') \subseteq \mathbb{Q}[Y]^0 \otimes_{\mathbb{Z}} (\kappa')^{\times}$, where $\mathbb{Q}[Y]^0 \subseteq \mathbb{Q}[Y]$ denotes the $\mathbb{Q}$-vector space of divisors of degree zero. 
    Moreover, we also have that $\operatorname{Im}(\delta') = \mathbb{Q}[Y]^0 \otimes_{\mathbb{Z}} (\kappa')^{\times}$, because $\delta'$ fits into the commutative diagram
    \[
        \begin{tikzcd}
            {H^{2,2}_{\mathcal{M}}((X \setminus Y)_{\kappa'})} \arrow[r, "\delta'"]                                                       & {\mathbb{Q}[Y]^0 \otimes_{\mathbb{Z}} (\kappa')^{\times}}                                                                                                                                                  \\
            {H^{1,1}_{\mathcal{M}}((X \setminus Y)_{\kappa'}) \otimes H^{1,1}_{\mathcal{M}}(\operatorname{Spec}(\kappa'))} \arrow[u, "\cup"] \arrow[r, "\sim"] & \left( \mathcal{O}^{\times}((X \setminus Y)_{\kappa'}) \otimes_{\mathbb{Z}} \mathbb{Q} \right) \otimes_{\mathbb{Z}} (\kappa')^{\times} \arrow[u, "\operatorname{div} \otimes \operatorname{Id}_{\kappa'}"]
        \end{tikzcd}
    \]
    and the divisor map 
    $\left( \mathcal{O}^{\times}((X \setminus Y)_{\kappa'}) \otimes_{\mathbb{Z}} \mathbb{Q} \right) \otimes_{\mathbb{Z}} (\kappa')^{\times} \xrightarrow{\operatorname{div} \otimes \operatorname{Id}_{\kappa'}} \mathbb{Q}[Y]^0 \otimes_{\mathbb{Z}} (\kappa')^{\times}$
    is surjective. This follows from the fact that we are taking rational coefficients, together with the assumption that there exists a point $y_0 \in Y$ such that $[\kappa(y_0) \colon \kappa] \, y - [\kappa(y) \colon \kappa] y_0 \in (J_X)_{\text{tors}}$ for every $y \in Y$.
    
    We have shown that
    \[
        H^{2,2}_{\mathcal{M}}((X \setminus Y)_{\kappa'}) \cong H^{2,2}_{\mathcal{M}}(X_{\kappa'}) \oplus \{ \mathcal{O}^{\times}((X \setminus Y)_{\kappa'}), (\kappa')^{\times} \}
    \]
    where $\{ \mathcal{O}^{\times}((X \setminus Y)_{\kappa'}), (\kappa')^{\times} \} \subseteq H^{2,2}_{\mathcal{M}}((X \setminus Y)_{\kappa'})$ denotes the subspace of symbols $\{f,c\} = \{ f \} \cup \{ c \}$ where $f \in \mathcal{O}^{\times}((X \setminus Y)_{\kappa'})$ and $c \in (\kappa')^{\times}$ is a constant. To conclude we can use Galois descent for motivic cohomology (see \cite[Theorem~14.3.4]{Cisinski_Deglise_2019}) to get an isomorphism
    \begin{equation} \label{eq:decomposition_motivic_cohomology}
        H^{2,2}_{\mathcal{M}}(X \setminus Y) \cong H^{2,2}_{\mathcal{M}}(X) \oplus \psi_{\ast}(\{ \mathcal{O}^{\times}((X \setminus Y)_{\kappa'}), (\kappa')^{\times} \})
    \end{equation}
    where $\psi \colon (X \setminus Y)_{\kappa'} \to X \setminus Y$ denotes the Galois covering induced by base change. Then the retraction 
    \[
        H^{2,2}_{\mathcal{M}}(X \setminus Y) \twoheadrightarrow H^{2,2}_{\mathcal{M}}(X)
    \]
    is simply given by the projection onto the first factor in the decomposition \eqref{eq:decomposition_motivic_cohomology}.
\end{proof}

We can now use the retraction $H^{2,2}_{\mathcal{M}}(X \setminus Y) \twoheadrightarrow H^{2,2}_{\mathcal{M}}(X)$ given by \cref{prop:motivic_cohomology_decomposition} to get a map
    \[
        \mathcal{O}^{\times}(X \setminus Y)^{\otimes 2} \otimes_{\mathbb{Z}} \mathbb{Q} \xrightarrow{\sim} H^{1,1}_{\mathcal{M}}(X \setminus Y) \xrightarrow{\cup} H^{2,2}_{\mathcal{M}}(X \setminus Y) \twoheadrightarrow H^{2,2}_{\mathcal{M}}(X)
    \]
which can be used to construct elements in motivic cohomology. 
This is a generalisation of ``Bloch's trick'' (see \cref{eq:bloch_trick}) that we will use in \cref{sec:defining_the_polynomial} to construct elements in the motivic cohomology group $H^{2,2}_{\mathcal{M}}(E)$ associated to an elliptic curve $E$ defined over $\mathbb{Q}$ which has potential complex multiplication.

\subsection{Deligne-Beilinson cohomology of curves over the reals} \label{sec:Deligne_cohomology}

In this section we recall the basic facts about Deligne-Beilinson cohomology groups of a smooth algebraic curve $X$ defined over $\mathbb{R}$. We will not discuss the general theory of Deligne-Beilinson cohomology. Suffice to say that the Deligne-Beilinson cohomology groups, denoted by $H^{i,j}_{\mathcal{D}}$, form a twisted Poincar{\'e} duality theory in the sense of \cite[\S 8]{Jannsen_1990}, which can be defined as: 
\begin{itemize}
    \item the hypercohomology of a suitable complex of sheaves (see \cite[Definition~2.9]{ev88});
    \item the sheaf cohomology of a suitable resolution of the previous complex (see \cite[Definition~5.50]{Burgos_Kramer_Kuhn_2007});
    \item the extension groups in the category of mixed Hodge structures (see \cite{Beilinson_1986});
    \item the cohomology induced by a motivic spectrum (see \cite[Section~6]{Bunke_Nikolaus_Tamme_2018}).
\end{itemize}
Deligne cohomology groups are also the target of Beilinson's \textit{regulator maps}
\[
    r_X^{\infty} \colon H^{i,j}_{\mathcal{M}}(X) \to H^{i,j}_{\mathcal{D}}(X)
\]
which can be constructed in many different ways, according to the chosen definition of motivic cohomology and Deligne cohomology (see \cite[\S~2]{be84}, \cite[\S~8]{Everest_Ward_1999}, \cite[\S~7]{Bunke_Nikolaus_Tamme_2018},\cite[Appendix]{Brunault_Zudilin_2020}).

We will only need the groups $H^{1,1}_{\mathcal{D}}(X)$ and $H^{2,2}_{\mathcal{D}}(X)$ for a smooth algebraic curve $X$ defined over $\mathbb{R}$. 
Hence we will only recall the definition of $H^{n,n}_{\mathcal{D}}(X)$ for a smooth variety $X$ defined over $\mathbb{R}$ or $\mathbb{C}$, following \cite[\S 7.3]{Nekovar_1994} (which is a special case of \cite[Definition~5.50]{Burgos_Kramer_Kuhn_2007}).

\begin{notation}
    We need to introduce the following notation:
    \begin{itemize}
        \item an analytic space $Y$ over $\mathbb{R}$ can be seen as a pair $(X,F_{\infty})$ where $X$ is a complex analytic space and $F_{\infty} \colon X \to X$ is an anti-holomorphic involution (see \cite[Teorema~14]{Tognoli_1967}). Moreover, a sheaf $\mathcal{S}$ on $Y$ can also be seen as a pair $(\mathcal{T},\sigma)$ where $\mathcal{T}$ is a sheaf on $X$ and $\sigma \colon F_{\infty}^{\ast}(\mathcal{T}) \to \mathcal{T}$ is an isomorphism whose inverse is $F_{\infty}^{\ast}(\sigma)$;
        \item for every algebraic variety $X$ over $\mathbb{C}$ we denote by $X(\mathbb{C})$ the usual complex analytification, given by the set of complex points endowed with the complex analytic topology. If $Y$ is an algebraic variety over $\mathbb{R}$ we denote by $Y^{\text{an}}$ the real analytic space $(Y_{\mathbb{C}}(\mathbb{C}),F_{\infty})$ where $F_{\infty}$ is complex conjugation (on points);
        \item for every subgroup $A \subseteq \mathbb{C}$ and every $j \in \mathbb{Z}$ we denote by $A(j) := (2 \pi i)^j \, A \subseteq \mathbb{C}$ and by $\pi_j \colon \mathbb{C} \to \mathbb{R}(j)$ the projection map given by $\pi_j(z) := (z + (-1)^j \, \overline{z})/2$. 
        If $X$ is a complex analytic space we denote by $\underline{A(j)}$ the constant sheaf with value $A(j)$, and if $Y = (X,F_{\infty})$ is a real analytic space we denote by $A(j)$ the pair $(\underline{A(j)},\overline{(\,)})$ where
        $\overline{(\,)} \colon F_{\infty}^{\ast}(\underline{A(j)}) = \underline{A(j)} \to \underline{A(j)}$
        denotes complex conjugation (on coefficients);
        \item for every smooth complex analytic space $X$ we denote by $\mathcal{A}^{\bullet,j}(X)$ the complex of smooth differential forms with values in $\mathbb{R}(j)$. If $Y$ is a smooth real analytic space given by the pair $(X,F_{\infty})$ we write $\mathcal{A}^{\bullet,j}(Y) := \mathcal{A}^{\bullet,j}(X)^{\overline{F_{\infty}^{\ast}}}$ where $\overline{(\,)}$ denotes again the action of complex conjugation on the coefficients of the differential forms. If $X$ is an algebraic variety over $\mathbb{C}$ (respectively, over $\mathbb{R}$) we write $\mathcal{A}^{\bullet,j}(X) := \mathcal{A}^{\bullet,j}(X(\mathbb{C}))$ (resp. $\mathcal{A}^{\bullet,j}(X) := \mathcal{A}^{\bullet,j}(X^{\text{an}})$);
        \item a \textit{good compactification} of a morphism $f \colon X \to Y$ of schemes (or analytic spaces) is a factorisation $f = p \circ j$ where $j \colon X \hookrightarrow Z$ is an open immersion, $p \colon Z \to Y$ is proper and $Z \setminus j(X)$ is a divisor with normal crossings. Moreover, if $f \colon X \to Y$ is smooth we assume that $p \colon Z \to Y$ is also smooth. When $Y = \operatorname{Spec}(\kappa)$ and $\kappa$ is a field of characteristic zero, we always have a good compactification, and any two good compactifications are dominated by a third one (see \cite[\S 3.2.II]{Deligne_1971}). When $X$ is a smooth curve over a field, then a good compactification is simply a smooth, proper curve $\overline{X}$ with an open immersion $j \colon X \hookrightarrow \overline{X}$ such that $\overline{X} \setminus j(X)$ is finite;
        \item if $\iota \colon D \hookrightarrow Z$ is a divisor with normal crossings on $Z$, and $j \colon Z \setminus D \hookrightarrow Z$ is the complementary open immersion, we denote by $\Omega^{\bullet}_{Z}\langle D \rangle \subseteq j_{\ast}(\Omega^{\bullet}_{Z \setminus D})$ the complex of sheaves of differential forms with logarithmic singularities along $D$ (see \cite[\href{https://stacks.math.columbia.edu/tag/0FUA}{Definition 0FUA}]{stacks-project}). This makes sense for schemes and also for analytic spaces. The global sections $\Omega^{\bullet}_{Z}\langle D \rangle(\overline{X}) \subseteq \Omega^{\bullet}_{Z \setminus D}(Z \setminus D)$ can be interpreted as (algebraic, smooth or holomorphic) differential forms on $Z \setminus D$ which have at worst logarithmic singularities ``at infinity'';
        \item for every smooth variety $X$ defined over $\mathbb{C}$ and any good compactification $X \hookrightarrow \overline{X}$ we define the complex
        \[
            \mathcal{F}^{\bullet}(X \hookrightarrow \overline{X}) := \Omega^{\bullet}_{\overline{X}(\mathbb{C})}\langle (\overline{X} \setminus X)(\mathbb{C})\rangle(\overline{X}(\mathbb{C}))
        \]
        which, up to quasi-isomorphism, is independent from the choice of a good compactification (see \cite[Theorem~5.46]{Burgos_Kramer_Kuhn_2007}). 
        For this reason, we will usually abuse notation and write $\mathcal{F}^{\bullet}(X) := \mathcal{F}^{\bullet}(X \hookrightarrow \overline{X})$;
        \item if $X$ is a smooth variety defined over $\mathbb{R}$ and $X \hookrightarrow \overline{X}$ is a good compactification we define the complex 
        \[
            \mathcal{F}^{\bullet}(X \hookrightarrow \overline{X}) := \mathcal{F}^{\bullet}\left( X_{\mathbb{C}} \hookrightarrow \overline{X}_{\mathbb{C}} \right)^{\overline{F_{\infty}^{\ast}}}
        \]
        and we will again abuse notation, denoting it by $\mathcal{F}^{\bullet}(X)$.
    \end{itemize} 
\end{notation}

\begin{definition} \label{def:simple_description_deligne_cohomology}
    Let $X$ be a smooth algebraic variety defined over $\mathbb{R}$ or $\mathbb{C}$. Then we define the Deligne cohomology groups
    \[
        H^{n,n}_{\mathcal{D}}(X) := \frac{\{ (\omega,\eta) \in \mathcal{A}^{n-1,n-1}(X) \oplus \mathcal{F}^n(X \hookrightarrow \overline{X}) \colon d(\omega) = \pi_{n-1}(\eta) \}}{d(\mathcal{A}^{n-2,n-1}(X))}
    \]
    where $X \hookrightarrow \overline{X}$ denotes any good compactification.
\end{definition}

\begin{remark} \label{rmk:cup_product_and_regulator_Deligne}
    We have an explicit description (see \cite[\S~3,10]{Everest_Ward_1999}) of the cup product
    \begin{align*}
        H^{n,n}_{\mathcal{D}}(X) \otimes H^{m,m}_{\mathcal{D}}(X) &\to H^{n+m,n+m}_{\mathcal{D}}(X) \\
        [(\omega_1,\eta_1)] \otimes [(\omega_2,\eta_2)] &\mapsto [(\omega_1 \wedge \pi_m(\eta_2) + (-1)^m \pi_n(\eta_1) \wedge \omega_2,\eta_1 \wedge \eta_2)]
    \end{align*}
    and of Beilinson's regulator map
    \begin{align*}
        r^{\infty}_X \colon \mathcal{O}^{\times}(X) \otimes_{\mathbb{Z}} \mathbb{Q} \cong H^{1,1}_{\mathcal{M}}(X) &\to H^{1,1}_{\mathcal{D}}(X) \\
        f \otimes 1 &\mapsto [(\log\lvert f \rvert,d\log(f))]
    \end{align*}
    which gives us the equality $r_X^{\infty}(\{f,g\}) = [(\log\lvert f \rvert d\operatorname{arg}(g) - \log\lvert g \rvert d\operatorname{arg}(f),0)]$.
\end{remark}

\begin{remark} \label{rmk:pairings}
    For every $n$-dimensional smooth algebraic variety $X$ over $\mathbb{R}$ or $\mathbb{C}$ we have an integration pairing
    \begin{align*}
        \langle \,, \, \rangle \colon  \mathcal{F}^j(X) \otimes H_j(X(\mathbb{C});\mathbb{R}) &\to \mathbb{C} \\
        \omega \otimes \gamma &\mapsto \langle \omega, \gamma \rangle := \int_{\gamma} \omega
    \end{align*}
    between differential forms and singular homology classes.
    If $X$ is proper then there is another pairing
        \begin{align*}
            \mathcal{A}^{n,j}(X) \otimes \mathcal{A}^{n,j}(X) &\to \mathbb{R} \\
            \alpha \otimes \beta &\mapsto \frac{1}{(\pi i)^j} \, \int_{X_{\mathbb{C}}(\mathbb{C})} \alpha \wedge \beta
        \end{align*}
    between differential forms, which is related to the first one by Poincar{\'e} duality (see \cite[\S~A.2.5]{Bost_1992}).
\end{remark}

Let now $Y$ be a smooth curve over $\mathbb{C}$, let $c \in \mathbb{C}^{\times}$ and let $f \in \mathcal{O}(Y)^{\times}$.
We can use the explicit descriptions provided by \cref{rmk:cup_product_and_regulator_Deligne} to compute the pairing of the regulator of the symbol $\{c,f\} \in H^{2,2}_{\mathcal{M}}(Y)$ with a homology class $\mathfrak{c} \in H_1(Y;\mathbb{Z})$.
To make this precise, let us recall
some elements from the theory of Riemann surfaces, following \cite[Appendix~A]{Bost_1992}.

\begin{remark} \label{rmk:riemann_surfaces}
    Let $X$ be a complex compact Riemann surface of genus $g$. Then the first singular homology group $H_1(X;\mathbb{Z})$ supports an intersection pairing $\# \colon H_1(X;\mathbb{Z})^{\otimes 2} \to \mathbb{Z}$ which is bilinear and anti-symmetric. 
    Moreover, $H_1(X;\mathbb{Z}) \cong \mathbb{Z}^{2 g}$, where $g \in \mathbb{N}$ denotes the genus of $X$, and there exists a $\mathbb{Z}$-basis $\{\alpha_i,\beta_j\}_{i,j=1,\dots,g} \subseteq H_1(X;\mathbb{Z})$ which is \textit{symplectic}, \textit{i.e.} for every $i, j \in \{1,\dots,g\}$ we have that 
    \[
        \alpha_i \# \alpha_i = \beta_j \# \beta_j = 0 \qquad \text{and} \qquad \alpha_i \# \beta_j = \delta_{i,j}
    \]
    where $\delta_{i,j} \in \{0,1\}$ denotes Kronecker's symbol (\textit{i.e.} $\delta_{i,j} = 1$ if and only if $i = j$).
    
    Now, let $S \subseteq X$ be a finite set of points and let $\iota \colon X \setminus S \hookrightarrow X$ denote the canonical inclusion. Then for every symplectic basis $\{\alpha_i,\beta_j\} \subseteq H_1(X;\mathbb{Z})$ and every point $x \in X \setminus S$ there exist  smooth loops $\{ a_i, b_j \colon [0,1] \to X \setminus S \}_{i,j = 1,\dots,g}$ such that:
    \begin{itemize}
        \item $a_i(0) = b_j(0) = a_i(1) = b_j(1) = x$ for every $i, j \in \{1,\dots,g\}$;
        \item $a_i( ]0,1[ ) \cap b_j(]0,1[) = \emptyset$ for every $i, j \in \{1,\dots,g\}$;
        \item $\restr{a_i}{[0,1[}$ and $\restr{b_j}{[0,1[}$ are injective for every $i,j \in \{1,\dots,g\}$;
        \item the vectors $\{a_i'(0), b_j'(0), a_i'(1), b_j'(1)\}_{i,j \in 1,\dots,g} \subseteq T_x(X)$ are pairwise non-collinear;
        \item the loops $\iota \circ a_i$ and $\iota \circ b_j$ are representatives of the homology classes $\alpha_i, \beta_j \in H_1(X;\mathbb{Z})$.
    \end{itemize}
    We will slightly abuse notation and denote by $\alpha_i, \beta_j \in H_1(X(\mathbb{C}) \setminus S; \mathbb{Z})$ the classes associated to the loops $a_i, b_j \colon [0,1] \to X \setminus S$.
    
    Now, observe that the loops $a_i, b_j$ correspond to a \textit{canonical dissection} $(\Delta,\varphi)$ of $X$ with $S \subseteq \varphi(\Delta^{\circ})$. More precisely, for every choice of $\{a_i, b_j\}$ as above there exists a polygon $\Delta \subseteq \mathbb{R}^2$ with $4g$ edges, an open $U \subseteq \mathbb{R}^2$ such that $\Delta \subseteq U$ and a surjective smooth map $\varphi \colon U \twoheadrightarrow X$ such that $\restr{\varphi}{\Delta^{\circ}}$ is a diffeomorphism onto $X \setminus C$ where
    \[
        C := \bigcup_{i} a_i([0,1]) \cup \bigcup_{j} b_j([0,1])
    \]
    is the union of all the loops given by $a_i$ and $b_j$. Each loop $a_i$ or $b_j$ corresponds to precisely two edges of $\Delta$ under $\varphi$, which are glued together with the same orientation (see \cite[Figure~23]{Bost_1992}).
    
    To conclude observe that for every $s \in S$ we can define a loop $c_s \colon [0,1] \to \Delta \setminus \varphi^{-1}(S) \to X \setminus S$, where the map $[0,1] \to \Delta \setminus \varphi^{-1}(S)$ is a small circle around $\varphi^{-1}(s)$ connected to one vertex of $\Delta$ by a straight line. 
    Let $\gamma_s \in H_1(X \setminus S;\mathbb{Z})$ be the singular cohomology class associated to $c_s$, which does not depend on the choice of the small circle $c_s$ if all the circles $\{c_s\}_{s \in S}$ are pairwise disjoint and oriented coherently.
    Then we have an exact sequence
    \begin{equation} \label{eq:homology_punctured_curve}
        \begin{aligned}
            0 \to \mathbb{Z} \to \mathbb{Z}^{S} &\to H_1(X(\mathbb{C}) \setminus S;\mathbb{Z}) \to H_1(X(\mathbb{C});\mathbb{Z}) \to 0 \\
            \{m_s\}_{s \in S} &\mapsto \sum_{s \in S} m_s \gamma_{s}
        \end{aligned}
    \end{equation}
    where the map $\mathbb{Z} \to \mathbb{Z}^{S}$ is the diagonal one. In particular, for every $s_0 \in S$ the set $\{\alpha_i, \beta_j, \gamma_s\}_{\substack{i,j = 1,\dots,g \\ s \in S \setminus \{s_0\}}}$ is a basis of $H_1(X \setminus S;\mathbb{Z})$. This can be easily shown using the Mayer-Vietoris exact sequence (see \cite[\S~4.6]{Spanier_1995}).
\end{remark}

Let us now use \cref{rmk:riemann_surfaces} to compute the pairing that we announced.

\begin{proposition} \label{prop:regulator_open_curve}
    Let $X$ be a smooth, proper algebraic curve over $\mathbb{C}$ of genus $g$, and let $S \subseteq X(\mathbb{C})$ be a finite set of points. Let $(\Delta,\varphi)$ be a canonical dissection of $X(\mathbb{C})$ such that $S \subseteq \varphi(\Delta^{\circ})$ and let $\alpha_i, \beta_j, \gamma_s$ be the homology classes associated to $(\Delta,\varphi)$. 
    Then we have that
    \begin{align}
        \langle r_{X \setminus S}^{\infty}(\{c,f\}), \alpha_i \rangle &= \langle r_{X \setminus S}^{\infty}(\{c,f\}), \beta_j \rangle = 0 \label{eq:regulator_from_compact}\\
        \langle r_{X \setminus S}^{\infty}(\{c,f\}), \gamma_s \rangle &= \log\lvert \partial_s(\{c,f\}) \rvert = \operatorname{ord}_s(f) \, \log\lvert c \rvert \label{eq:regulator_tame_symbol}
    \end{align}
    for every $c \in \mathbb{C}$, every $f \in \mathbb{C}(X)$ such that $S_f \subseteq S$, every $i, j \in \{1,\dots,g\}$ and every $s \in S$.
\end{proposition}
\begin{proof}
    The computation \eqref{eq:regulator_from_compact} follows from the fact that $\int_{X(\mathbb{C})} d(\log(f)) \wedge \alpha = 0$, whereas \eqref{eq:regulator_tame_symbol} is an application of \textit{Jensen's formula}, as explained in \cite[Page~25]{vi99}.
\end{proof}

To conclude this preliminary section let us introduce some notations concerning the cohomology of elliptic curves defined over the reals.
\begin{notation}
    Let $E$ be an elliptic curve defined over $\mathbb{R}$. We introduce the following notation:
    \begin{itemize}
        \item $E(\mathbb{R})^0 \subseteq E(\mathbb{R})$ denotes the connected component of the identity;
        \item $\omega_E \in \mathcal{F}^1(E)$ the unique differential form such that $\int_{E(\mathbb{R})^{0}} \omega_E = 1$. Observe that $\omega_E \in H^1(E^{\text{an}};\mathbb{Q}(1))$ because $H_1(E^{\text{an}};\mathbb{Q})$ is generated by the homology class of $E(\mathbb{R})^0$;
        \item $H_1(E(\mathbb{C});\mathbb{Q})^{-} \subseteq H_1(E(\mathbb{C});\mathbb{Q})$ denotes the subspace of homology classes which are anti-invariant by complex conjugation;
        \item $\gamma_E \in H_1(E(\mathbb{C});\mathbb{Q})^{-}$ denotes the Poincar{\'e} dual of $\omega_E$.
    \end{itemize}
\end{notation}

\section{Constructing the polynomials}
\label{sec:defining_the_polynomial}

The aim of this section is to associate to every elliptic curve $E$ defined over $\mathbb{Q}$ which has potential complex multiplication by the ring of integers $\mathcal{O}_K$ of an imaginary quadratic field $K$ the polynomial $P \in \mathbb{Z}[x,y]$ appearing in \cref{thm:main_theorem}. 
To do so, we will study pairs of functions $f,g \colon E \to \mathbb{P}^1$ defined in \cite[Theorem~4.10]{dw88} and \cite[Page~384]{ro87}, and we will prove that $\mathbb{Q}(E) = \mathbb{Q}(f,g)$. Hence if we take $P \in \mathbb{Z}[x,y]$ to be the minimal polynomial of $f$ and $g$ we will see immediately that $V_P$ is birational to $E$, which was one of the conditions outlined in the statement of \cref{thm:main_theorem}. 
The Mahler measure $m(P)$ of $P$ will be related to $L'(E,0)$ in \cref{sec:mahler_measure}, which will give a complete proof of \cref{thm:main_theorem}.

The pairs of functions $f,g \colon E \to \mathbb{P}^1$ that we are looking for have the property that $S_{f,g} \subseteq E(\overline{\mathbb{Q}})_{\text{tors}}$. Hence $n_{f,g} \in \mathbb{Z}_{\geq 1}$ and we can define, following Bloch (see \cite[Proposition~10.1.1]{bl00}), a motivic cohomology class
\begin{equation} \label{eq:bloch_trick}
\eta_{f,g} := n_{f,g} \, \{ f,g \} + \sum_{x \in S_{f,g} \setminus \{ 0 \}} \{ \partial_x(\{f,g\}),\varphi_{f,g}^{(x)} \} \in H^{2,2}_{\mathcal{M}}(E)
\end{equation}
which has the remarkable property that $\langle r^{\infty}_{E}(\eta_{f,g}), \gamma_E \rangle/L'(E,0) \in \mathbb{Q} $ (see \cref{sec:preliminaries} for all the relevant notation). This property is made explicit by the following result, which is due to Rohrlich.

\begin{theorem}[see {\cite{ro87}}] \label{thm:rohrlich_general_theorem}
    Let $E$ be an elliptic curve defined over $\mathbb{Q}$ having potential complex multiplication by the ring of integers $\mathcal{O}_K$ of an imaginary quadratic field $K$. Let moreover $f, g \colon E \to \mathbb{P}^1$ be two functions such that $S_{f,g} \subseteq E(\overline{\mathbb{Q}})_{\text{tors}}$. Then we have that 
    \[
        \langle r^{\infty}_E(\eta_{f,g}), \gamma_E \rangle = \mathcal{R}(\operatorname{div}(f) \diamondsuit \operatorname{div}(g)) \, L'(E,0)
    \] 
    and $n_{f,g} \, \mathcal{R}(\operatorname{div}(f) \diamondsuit \operatorname{div}(g)) \in \mathbb{Z}$ (see \cref{sec:Deligne_cohomology} and \cref{sec:CM_elliptic_curves} for the relevant definitions).
\end{theorem}

Hence to prove Beilinson's conjectures for the special value $L^{\ast}(E,0) = L'(E,0)$ one has to show that for every CM elliptic curve $E$ defined over $\mathbb{Q}$ we can find a pair of functions $f,g \colon E \to \mathbb{P}^1$ such that $\mathcal{R}(\operatorname{div}(f) \diamondsuit \operatorname{div}(g)) \neq 0$. This happens for many pairs of functions, as we will explain in \cref{sec:deninger_wingberg} and \cref{sec:rohrlich_polynomials}.

Before doing that, let us make some remarks concerning the construction \eqref{eq:bloch_trick}, which sometimes goes under the name of \textit{``Bloch's trick''}.
\begin{remark}
    It is easy to see that $\eta$ is bilinear, alternating and invariant by scaling, \textit{i.e.}
    \[
        \eta_{fg,h} = \eta_{f,h} + \eta_{g,h}, \quad \text{and} \quad \eta_{f,g} = -\eta_{g,f} \quad \text{and} \quad \eta_{c,f} = 0
    \]
    for every $f,g \in \mathbb{Q}(E)$ and $c \in \mathbb{Q}$. This shows that we have an alternating, bilinear pairing 
    \begin{equation} \label{eq:map_on_divisors}
        \begin{aligned}
            [\,,\,]_{\mathcal{M}} \colon \bigwedge\nolimits^2 \mathbb{Q}[E(\overline{\mathbb{Q}})_{\text{tors}}]^{0,\operatorname{Gal}(\overline{\mathbb{Q}}/\mathbb{Q})} &\to H^2_{\mathcal{M}}(E ; \mathbb{Q}(2)) \\
            D_1 \wedge D_2 &\mapsto \eta_{f_1,f_2}
        \end{aligned}
    \end{equation}
    where $f_1,f_2 \colon E \to \mathbb{P}^1$ are any two functions such that $\operatorname{div}(f_j) = \operatorname{ord}(D_j) \, D_j$ (see \cref{not:orders_and_divisors}). 
    
    Observe finally that for every $D_1, D_2 \in \mathbb{Q}[E(\overline{\mathbb{Q}})_{\text{tors}}]^{0,\operatorname{Gal}(\overline{\mathbb{Q}}/\mathbb{Q})}$ we have that 
    \[
        [D_1,D_2]_{\mathcal{M}} = n_1 \, n_2 \, [D_1,D_2]_{\mathcal{A}}
    \]
    where $[\,,\,]_{\mathcal{A}}$ is the pairing defined in \cite[Theorem~5.1]{dw88}.
\end{remark}

\subsection{Models of CM elliptic curves (according to Deninger and Wingberg)}
\label{sec:deninger_wingberg}

The aim of this section is to construct the first pair of functions $f,g \in \mathbb{Q}(E)$ of the kind described at the beginning of this section. Let us start with the following result, which is due to Deninger and Wingberg.

\begin{lemma}[see {\cite[Theorem~4.10]{dw88}}] \label{lem:deninger_wingberg_construction}
    Let $E$ be an elliptic curve defined over $\mathbb{Q}$ having potential complex multiplication by the ring of integers $\mathcal{O}_K$ of an imaginary quadratic field $K$. 
    Then there exist two functions $f,g \colon E \to \mathbb{P}^1$ such that
    \begin{align*}
        \operatorname{div}(f) &= \sum_{x \in E[\nu_E](\overline{\mathbb{Q}}) \setminus \{0\}} \left( (x) - (0) \right) \\
        \operatorname{div}(g) &= c_g \, \sum_{y \in E[\nu_E](\overline{\mathbb{Q}})/\mu_K} \left( ([\widetilde{\chi_E}(\overline{x})](x)) - (0) \right)
    \end{align*}
    where $c_g \in \{1,2\}$ denotes the order of the point $\sum_{y \in E[\nu_E](\overline{\mathbb{Q}})/\mu_K} [\widetilde{\chi_E}(\overline{x})](x) \in E[2](\mathbb{Q})$.
    Moreover we have that 
    \begin{equation} \label{eq:deninger_wingberg_rational_factor}
        \mathcal{R}(\operatorname{div}(f) \diamondsuit \operatorname{div}(g)) = \frac{c_g \, \mathfrak{f}_E}{\lvert \operatorname{disc}(K/\mathbb{Q}) \rvert} = c_g \, \operatorname{N}_{K/\mathbb{Q}}(\mathfrak{f}_{\psi_{E_K}}) \in \mathbb{Z} \setminus \{0\}
    \end{equation}
    where $E_K$ is the base change of $E$ over $K$ and $\psi_{E_K}$ denotes the Hecke character defined in \cref{sec:CM_elliptic_curves}. 
\end{lemma}
\begin{proof}
    The two divisors
    \[
        \sum_{x \in E[\nu_E](\overline{\mathbb{Q}}) \setminus \{0\}} \left( (x) - (0) \right) \qquad \text{and} \qquad 2 \, \sum_{y \in E[\nu_E](\overline{\mathbb{Q}})/\mu_K} \left( ([\widetilde{\chi_E}(\overline{x})](x)) - (0) \right)
    \]
    are elements of $\mathbb{Q}[E(\overline{\mathbb{Q}})_{\text{tors}}]^{0,\operatorname{Gal}(\overline{\mathbb{Q}}/\mathbb{Q})}$, as it is clear from the explicit description of the Galois action on torsion points (see \cite[Section~4]{dw88} and \cref{sec:CM_elliptic_curves}). Moreover, the fact that $E[\nu_E](\overline{\mathbb{Q}})$ is a group implies that
    \[
        \sum_{x \in E[\nu_E](\overline{\mathbb{Q}}) \setminus \{0\}} x = \begin{cases} 0, \ \text{if} \ 2 \nmid \operatorname{N}_{K/\mathbb{Q}}(\nu_E) \\ \sum_{x \in E[2](\overline{\mathbb{Q}}) \setminus \{0\}} x = 0, \ \text{otherwise} \end{cases}
    \]
    which follows from the fact that $E[2](\overline{\mathbb{Q}}) \cong (\mathbb{Z}/2 \mathbb{Z})^2$. For similar reasons we have that
    \[
        \sum_{y \in E[\nu_E](\overline{\mathbb{Q}})/\mu_K} [\widetilde{\chi_E}(\overline{x})](x) \in E[2](\mathbb{Q})
    \]
    which implies that we can find two functions $f,g \colon E \to \mathbb{P}^1$ as in the statement of the theorem. 
    Now the identity \eqref{eq:deninger_wingberg_rational_factor} follows from the computations carried out in \cite[Section~4]{dw88}, after having observed that the regulator used by Rohrlich is twice the regulator used by Deninger and Wingberg (see \cite[Page~371]{ro87} and \cite[Equation~1.8]{dw88} for a comparison) and that $\operatorname{div}(g)$ is twice the divisor $\beta$ which appears in \cite[Theorem~4.10]{dw88}.
\end{proof}

\begin{remark}
    It would in principle be possible to prove the identity \eqref{eq:deninger_wingberg_rational_factor} using directly the definition of $\mathcal{R}$ (see \cref{not:rational_factor_function}). However this seems difficult, given the complexity of the divisors involved in \cref{lem:deninger_wingberg_construction}.
\end{remark}

We will now use an idea due to Brunault (see \cite[Lemma~3.3]{br16}) to prove that $\mathbb{Q}(E) = \mathbb{Q}(f,g)$.

\begin{lemma} \label{lem:brunault_generating_functions}
    Let $E$ be an elliptic curve defined over a field $\kappa$. For every $P \in E(\overline{\kappa})_{\text{tors}}$ let $O_P := \operatorname{Gal}(\overline{\kappa}/\kappa) \, P$ and let $f_P \in \kappa(E)$ be any function such that
    \[
        \operatorname{div}(f_P) = c_P \, \sum_{x \in O_P} \left( (x) - (0) \right)
    \]
    where $c_P \in \mathbb{Z}_{\geq 1}$ is the order of the point $\sum_{x \in O_P} x \in E(\kappa)_{\text{tors}}$. 
    
    Then we have that: 
    \begin{enumerate}
        \item the extension $\kappa(f_P) \subset \kappa(E)$ contains no proper sub-extensions;
        \item if $\kappa(f_P) = \kappa(f_Q)$ for some points $P, Q \in E(\overline{\kappa})_{\text{tors}}$ and $\operatorname{char}(\kappa) = 0$ then $\lvert O_P \rvert = \lvert O_Q \rvert$.
    \end{enumerate}
\end{lemma}
\begin{proof}
    Consider a sub-extension $\kappa(f_P) \subseteq F \subseteq \kappa(E)$. Two possibilities can occur:
    \begin{itemize}
        \item $F = \kappa(g)$ for some function $g \in \kappa(E)$, which implies that $f_P = h \circ g$ for some $h \colon \mathbb{P}^1_{\kappa} \to \mathbb{P}^1_{\kappa}$. We can assume, up to applying two homographies $\mathbb{P}^1_{\kappa} \to \mathbb{P}^1_{\kappa}$, that $g(0) = \infty$ and that $h(0) = 0$. These homographies can be taken to be defined over $\kappa$ because $0 \in \mathbb{P}^1(\kappa)$ and $g(0) \in \mathbb{P}^1(\kappa)$. Then every zero of $g$ is a zero of $f_P$, and the converse also applies because $g$ is not constant (hence it has some zero $x \in O_P$) and defined over $\kappa$ (hence all the points $y \in O_x = O_P$ are zeros of $g$). Moreover $h(\infty) = h(g(0)) = f_P(0) = \infty$, which implies that $0$ is the unique pole of $g$ (since $0$ is the unique pole of $f_P$). This implies that
        \[
            \operatorname{div}(g) = d \, \sum_{x \in \operatorname{Gal}(\overline{\kappa}/\kappa) \, P} \left( (x) - (0) \right)
        \]
        for some $d \in \mathbb{Z}_{\geq 1}$. But then $c_P \mid d$ (since $c_P$ is the order of $\sum_{x \in O_P} x \in E(\kappa)_{\text{tors}}$) and thus $c_P = d$ (because $f_P = h \circ g$). Hence $g = \alpha \, f_P$ for some $\alpha \in \kappa^{\times}$, which implies that $F = \kappa(f_P)$.
        
        \item there is an isogeny $\varphi \colon E \twoheadrightarrow E'$ which induces an embedding $\varphi^{\ast} \colon \kappa(E') \hookrightarrow \kappa(E)$ and we have that $F = \varphi^{\ast}(\kappa(E'))$. This implies that $f_P = g \circ \varphi$ for some function $g \in \kappa(E')$, which in turn implies that $f_P(x) = \infty$ for every $x \in \ker(\varphi)$. Hence $\varphi$ is an isomorphism (because $0$ is the unique pole of $f_P$) and thus $F = \kappa(E)$.
    \end{itemize}
    This shows that $\kappa(f_P) \subset \kappa(E)$ contains no proper sub-extensions.
    
    Now suppose that $\kappa(f_P) = \kappa(f_Q)$ for some points $P, Q \in E(\overline{\kappa})_{\text{tors}}$. Then we have that 
    \[
    c_P \, \lvert O_P \rvert = [\kappa(E) \colon \kappa(f_P)] = [\kappa(E) \colon \kappa(f_Q)] = c_Q \, \lvert O_Q \rvert
    \] 
    (see \cite[Proposition~8.4]{fu89}) and that
    \[
        f_P = \frac{a f_Q + b}{c f_Q + d} \qquad \text{for some} \ \begin{pmatrix} a & b \\ c & d \end{pmatrix} \in \operatorname{GL}_2(\kappa)
    \]
    and since both $f_P$ and $f_Q$ have $0$ as their unique pole we must have that $c = 0$. Hence we get
    \begin{equation} \label{eq:riemann_hurwitz}
        \lvert O_P \rvert \, c_P + 1 \geq \sum_{x \in O_P} (c_P - 1) \, [\kappa(x) \colon \kappa] + \sum_{x \in O_Q} (c_Q - 1) \, [\kappa(x) \colon \kappa]
    \end{equation}
    applying the Riemann-Hurwitz formula (see \cite[\href{https://stacks.math.columbia.edu/tag/0C1B}{Section~0C1B}]{stacks-project}) for the covering $f_P \colon E \to \mathbb{P}^1$.
    This implies that $c_P = c_Q = 1$ and thus that $\lvert O_P \rvert = \lvert O_Q \rvert$.
\end{proof}

\begin{theorem} \label{thm:minimal_polynomial_deninger_wingberg}
    Let $E$ be an elliptic curve defined over $\mathbb{Q}$ having potential complex multiplication by the ring of integers $\mathcal{O}_K$ of an imaginary quadratic field $K$. 
    Let moreover $f,g \in \mathbb{Q}(E)$ be as in \cref{lem:deninger_wingberg_construction}.
    
    Then we have that $\mathbb{Q}(E) = \mathbb{Q}(f,g)$ and $\deg_x(P) = \operatorname{N}_{K/\mathbb{Q}}(\nu_E) - 1$, where $P \in \mathbb{Z}[x,y]$ denotes any minimal polynomial for $f$ and $g$.
\end{theorem}
\begin{proof}
    We know that $[\mathbb{Q}(E) \colon \mathbb{Q}(f)] = \lvert E[\nu_E](\overline{\mathbb{Q}}) \setminus \{ 0 \} \rvert = \operatorname{N}_{K/\mathbb{Q}}(\nu_E) - 1$ (see \cite[Proposition~8.4]{fu89}), which implies that $\deg_x(P) = \operatorname{N}_{K/\mathbb{Q}}(\nu_E) - 1$. Moreover, $[\mathbb{Q}(E) \colon \mathbb{Q}(g)] < [\mathbb{Q}(E) \colon \mathbb{Q}(f)]$ because $\lvert E[\nu_E](\overline{\mathbb{Q}})/\mu_K \rvert < \lvert E[\nu_E](\overline{\mathbb{Q}}) \rvert$. 
    
    We also have that $\mathbb{Q}(g) = \mathbb{Q}(f_P)$ where $P = [\widetilde{\chi_E}(\overline{x_0})](x_0)$ for any $x_0 \in E[\nu_E](\overline{\mathbb{Q}})$. Indeed, we know that for every $x \in E[\nu_E](\overline{\mathbb{Q}})$ there exists $a \in (\mathcal{O}_K/\nu_E)^{\times}$ such that $x = [a^{-1}](x)$, because $E[\nu_E](\overline{\mathbb{Q}})$ is a free $(\mathcal{O}_K/\nu_E)$-module of dimension one (see \cite[Lemma~1]{pa89}). We can now use \cref{lem:galois_action_on_torsion} to see that
    \begin{align*}
        \operatorname{div}(g) &= c_g \, \sum_{y \in E[\nu_E](\overline{\mathbb{Q}})/\mu_K} \left( ([\widetilde{\chi_E}(\overline{x})](x)) - (0) \right) = \\
        &= c_g \, \sum_{a \in (\mathcal{O}_K/\nu_E)^{\times}/\mu_K} \left( ([\widetilde{\chi_E}(\overline{[a^{-1}](x_0)})]([a^{-1}](x_0))) - (0) \right) = \\
        &= c_g \, \sum_{a \in (\mathcal{O}_K/\nu_E)^{\times}/\mu_K} \left( ([a^{-1} \, \chi_E(a) \, \widetilde{\chi_E}(\overline{x})](x)) - (0) \right) = \\
        &= c_g \, \sum_{\sigma \in \operatorname{Gal}(K(E[\nu_E](\overline{\mathbb{Q}}))/K)} \left( (\sigma([\widetilde{\chi_E}(\overline{x})](x))) - (0) \right) = c_P \, \sum_{y \in O_P} \left( (y) - (0) \right) = \operatorname{div}(f_P)
    \end{align*}
    which implies that $g = \alpha \, f_P$ for some $\alpha \in \mathbb{Q}^{\times}$.
    
    Now to conclude that $\mathbb{Q}(E) = \mathbb{Q}(f,g)$ we can apply \cref{lem:brunault_generating_functions}, using the fact that $\mathbb{Q}(f) \neq \mathbb{Q}(f_P)$ since $[\mathbb{Q}(E) \colon \mathbb{Q}(f_P)] = [\mathbb{Q}(E) \colon \mathbb{Q}(g)] < [\mathbb{Q}(E) \colon \mathbb{Q}(f)]$. 
\end{proof}

\begin{remark}
    We know that $\deg_x(P) = \operatorname{N}_{K/\mathbb{Q}}(\nu_E) - 1$. 
    Computing $\deg_y(P)$ is harder, but it can be done if we know $\lvert E[\nu_E](\overline{\mathbb{Q}})/\mu_K \rvert$ (which depends on $\gcd(\operatorname{N}_{K/\mathbb{Q}}(\nu_E),\lvert \mu_K \rvert)$) and $\lvert S \rvert$, where
    \[
        S := \{ x \in E[\nu_E](\overline{\mathbb{Q}}) \mid \widetilde{\chi_E}(x) = 0 \} = \bigcup_{\substack{\alpha \mid \mathfrak{f}_{\psi_{E_K}} \\ \alpha \neq 1}} E\left[ \frac{\nu_E}{\alpha} \right](\overline{\mathbb{Q}}) =  \bigcup_{\substack{\alpha \mid \mathfrak{f}_{\psi_{E_K}} \\ \alpha \, \mathcal{O}_K \in \operatorname{Spec}(\mathcal{O}_K) }} E\left[ \frac{\nu_E}{\alpha} \right](\overline{\mathbb{Q}})
    \]
    which shows that $\lvert S \rvert$ can be computed using an inclusion-exclusion principle.
\end{remark}

\subsection{Models of CM elliptic curves (according to Rohrlich)}
\label{sec:rohrlich_polynomials}
Let us turn our attention to the pair(s) of functions $f,g$ constructed by Rohrlich.
We keep again using the notation introduced in \cref{sec:preliminaries}.

\begin{lemma}[see {\cite[Pages~384-386]{ro87}}] \label{lem:rohrlich_functions}
    Let $E$ be an elliptic curve defined over $\mathbb{Q}$ having potential complex multiplication by the ring of integers $\mathcal{O}_K$ of an imaginary quadratic field $K$. 
    Let $p \in \mathbb{N}$ be a prime such that $p \nmid \mathfrak{f}_E$ and $p \, \mathcal{O}_K$ is also prime. 
    Let moreover $c \in \mathbb{N}$ be an integer such that 
    \[
        \mathfrak{f}_{\psi_{E_K}} \, \mathfrak{b}_E^{-1} \mid c \, \mathcal{O}_K \mid \mathfrak{f}_{\psi_{E_K}}^m \qquad \text{for some} \qquad m \in \mathbb{N}
    \]
    where $E_K$ denotes the base change of $E$ to $K$ and all the other terms are defined in \cref{sec:CM_elliptic_curves}. 
    
    Then there exist two functions $f,g \colon E \to \mathbb{P}^1$ such that
    \begin{align*}
        \operatorname{div}(f) &= k_p \, \sum_{x \in O_p} \left( (x) - (0) \right) \\
        \operatorname{div}(g) &= k_c \, \sum_{y \in O_c} \left( (y) - (0) \right)
    \end{align*}
    where for every $m \in \mathbb{Z}_{\geq 1}$ we define $O_{m} := \operatorname{Gal}(\overline{\mathbb{Q}} / \mathbb{Q}) \, \theta_E(1/m) \subseteq E[m](\overline{\mathbb{Q}})$ and $k_m \in \mathbb{Z}_{\geq 1}$ to be the order of the torsion point $\sum_{x \in O_m} x \in E[m](\mathbb{Q})$. Finally, we have that
    \begin{equation} \label{eq:rational_factor_Rohrlich}
        \mathcal{R}\left( \operatorname{div}(f) \diamondsuit \operatorname{div}(g) \right) = - \frac{k_p \, k_c \, (1 + p^3)}{c \, p} \in \frac{1}{n_{f,g}} \, \mathbb{Z} \setminus \{ 0 \}.
    \end{equation}
\end{lemma}
\begin{proof}
    First of all, observe that such a number $c \in \mathbb{N}$ exists because $\mathfrak{b}_E^{-1} \mid \mathfrak{f}_{\psi_{E_K}}$, which follows from Deuring's formula (see \cref{prop:deuring_CM_conductor}) and the fact that $\operatorname{ord}_p(\mathfrak{f}_E) \neq 1$ for every prime $p \in \mathbb{N}$.
    Now observe that $\sum_{x \in O_m} \left( (x) - (0) \right) \in \mathbb{Q}[E(\overline{\mathbb{Q}})_{\text{tors}}]^{0,\operatorname{Gal}(\overline{\mathbb{Q}}/\mathbb{Q})}$ for every $m \in \mathbb{Z}_{\geq 1}$, which implies the existence of the pair $f,g \in \mathbb{Q}(E)$.
    
    Let us now turn to the proof of \eqref{eq:rational_factor_Rohrlich}. First of all, it is evident from the definition that $p \nmid c$, which implies that for every $x \in E[p](\overline{\mathbb{Q}})$ and $y \in E[c](\overline{\mathbb{Q}})$ we have that $\lvert \operatorname{Gal}(\overline{\mathbb{Q}}/\mathbb{Q}) \, (x,y) \rvert = \lvert \operatorname{Gal}(\overline{\mathbb{Q}}/\mathbb{Q}) \, (x - y) \rvert$. Moreover, for every $\sigma_1, \sigma_2 \in \operatorname{Gal}(\overline{\mathbb{Q}}/\mathbb{Q})$ there exists $\tau \in \operatorname{Gal}(\overline{\mathbb{Q}}/\mathbb{Q})$ such that
    \[
        \sigma_1\left( \theta_E\left( \frac{1}{p} \right) \right) = \tau \left( \theta_E\left( \frac{1}{p} \right) \right) \qquad \text{and} \qquad \sigma_2\left( \theta_E\left( \frac{1}{c} \right) \right) = \tau \left( \theta_E\left( \frac{1}{c} \right) \right)
    \]
    because $\operatorname{Gal}(\overline{\mathbb{Q}}/\mathbb{Q}) \, \theta_E(\alpha) = \operatorname{Gal}(\overline{\mathbb{Q}}/K) \, \theta_E(\alpha)$ for every $\alpha \in \mathbb{R}$ and $K(E[p](\overline{\mathbb{Q}})) \cap K(E[c](\overline{\mathbb{Q}})) = K$ (see \cref{prop:campagna_pengo}). 
    This implies that
    \begin{equation} \label{eq:rohrlich_proof_1}
        \begin{aligned}
        \mathcal{R}\left( \operatorname{div}(f) \diamondsuit \operatorname{div}(g) \right) = k_p \, k_c \, ( &\mathcal{R}(\theta_E(1/p) - \theta_E(1/c)) - \lvert O_p \rvert \, \mathcal{R}(\theta_E(1/c)) - \\ 
        &- \lvert O_c \rvert \, \mathcal{R}(\theta_E(-1/p)) + \lvert O_p \rvert \, \lvert O_c \rvert \, \mathcal{R}(0) )
    \end{aligned}
    \end{equation}
    because $\{ (0,0), (\theta_E(1/p),0), (0,\theta_E(1/c)), (\theta_E(1/p),\theta_E(1/c)) \}$ is a full set of representatives for the diagonal action of $\operatorname{Gal}(\overline{\mathbb{Q}}/\mathbb{Q})$ on $S_f \times S_g$. We have moreover that $\mathcal{R}(\theta_E(-1/p)) = \mathcal{R}(0) = 0$ and that $\lvert O_p \rvert = p^2 - 1$ (see \cite[Theorem~7.8(c)]{bc16}). Observe now that
    \begin{equation} \label{eq:rohrlich_proof_2}
        \mathcal{R}(\theta_E(1/c)) = \frac{1}{c} \, \prod_{\mathfrak{p} \mid c \, \mathfrak{b}_E} \left( 1 - \psi_{E_K}(\mathfrak{p}) \right) = \frac{1}{c}
    \end{equation}
    because $\theta_E(1/c) \in E(\mathbb{R})$, no prime ideal $\mathfrak{p} \mid c \, \mathfrak{b}_E$ is coprime to $\mathfrak{f}_{\psi_{E_K}}$ and $\operatorname{Ann}_{\mathcal{O}_K}(\theta_E(1/c)) = c \, \mathfrak{b}_E$. Finally, we have that
    \begin{equation} \label{eq:rohrlich_proof_3}
        \begin{aligned}
        \mathcal{R}(\theta_E(1/p) - \theta_E(1/c)) &= \mathcal{R}\left( - \theta_E\left( \frac{\psi_{E_K}((c - p) \, \mathcal{O}_K)}{c \, p} \right) \right) = \frac{-1}{c \, p} \, \prod_{\mathfrak{p} \mid (c \, p) \, \mathfrak{b}_E} \left( 1 - \psi_{E_K}(\mathfrak{p}) \right) = \\
        &= \frac{-(1 - \psi_{E_K}(p \, \mathcal{O}_K))}{c \, p} = - \frac{1 + p}{c \, p}
    \end{aligned}
    \end{equation}
    because $\operatorname{Ann}_{\mathcal{O}_K}(\theta_E(1/(c \, p))) = c \, p \, \mathfrak{b}_E$ and the only prime which divides $c \, p \, \mathfrak{b}_E$ and is coprime with $\mathfrak{f}_{\psi_{E_K}}$ is $p \, \mathcal{O}_K$, for which we have that $\psi_{E_K}(p \, \mathcal{O}_K) = \left( \frac{\operatorname{disc}(K/\mathbb{Q})}{p} \right) \, p = -p$. Putting together \eqref{eq:rohrlich_proof_1}, \eqref{eq:rohrlich_proof_2} and \eqref{eq:rohrlich_proof_3} we obtain \eqref{eq:rational_factor_Rohrlich}.
\end{proof}

\begin{remark}
    Observe that $k_m \in \{1,2,3,4,6\}$ for every $m \in \mathbb{Z}_{\geq 1}$, which follows from the complete characterisation of the possible rational torsion subgroups $E(\mathbb{Q})_{\text{tors}}$ associated to an elliptic curve $E$ defined over $\mathbb{Q}$ which has potential complex multiplication (see \cite{ol74}).
\end{remark}

\begin{remark}
    If we take $c \in \mathbb{Z}_{\geq 1}$ such that $\mathfrak{f}_{\psi_{E_K}} \mid c \, \mathcal{O}_K$ we know that $K(E[c](\overline{\mathbb{Q}}))$ coincides with the ray class field of $K$ relative to the modulus $c \, \mathcal{O}_K$ (see \cite[Lemma~3]{Coates_Wiles_1977}). Hence in this case we do not need to use \cref{prop:campagna_pengo} to prove \cref{lem:rohrlich_functions}.
\end{remark}

We can now prove the analogue of \cref{thm:minimal_polynomial_deninger_wingberg} for Rohrlich's functions.

\begin{theorem} \label{thm:minimal_polynomial_rohrlich}
Let $E$ be an elliptic curve defined over $\mathbb{Q}$ having potential complex multiplication by the ring of integers $\mathcal{O}_K$ of an imaginary quadratic field $K$. 
Let $p, c \in \mathbb{N}$ and $f, g \in \mathbb{Q}(E)$ be as in \cref{lem:rohrlich_functions} and assume that $\phi(c) > p^2 - 1$, where $\phi$ denotes Euler's totient function. Then $\mathbb{Q}(f,g) = \mathbb{Q}(E)$ and if $P \in \mathbb{Z}[x,y]$ denotes a minimal polynomial of $f$ and $g$ we have that $\deg_x(P) = \lvert O_p \rvert$ and $\deg_y(P) = \lvert O_c \rvert$.
\end{theorem}
\begin{proof}
    We see from \cref{lem:brunault_generating_functions} that either $\mathbb{Q}(E) = \mathbb{Q}(f,g)$ or $\mathbb{Q}(f) = \mathbb{Q}(g)$, and in this case we would have that $\lvert O_p \rvert = \lvert O_c \rvert$, but this is absurd. Indeed, $\lvert O_p \rvert = p^2 - 1$ and $\phi(c) < \lvert O_c \rvert$ (see \cite[Section~6.5]{bc16}). Then our hypothesis shows that $\lvert O_c \rvert > \lvert O_p \rvert$.
    The final part of the theorem follows simply from the fact that $\deg_x(P) = [\mathbb{Q}(E) \colon \mathbb{Q}(f)]$ and $\deg_y(P) = [\mathbb{Q}(E) \colon \mathbb{Q}(g)]$.
\end{proof}

\section{Computing the Mahler measure} \label{sec:mahler_measure}

The aim of this section is to complete the proof of \cref{thm:main_theorem}, taking as $P \in \mathbb{Z}[x,y]$ a slightly modified version of the polynomials that we defined in \cref{sec:defining_the_polynomial}. To do so observe that for every ring $R$ and every $n \in \mathbb{Z}_{\geq 1}$ the ring of Laurent polynomials $R[x_1^{\pm 1},\dots,x_n^{\pm 1}] = \Gamma(\mathbb{G}_{m,R}^n,\mathcal{O}_{\mathbb{G}_{m,R}^n})$ supports the action of the group
\[
    (R^{\times})^{n+1} \times \operatorname{Aut}(\mathbb{G}_{m}^n) \cong (R^{\times})^{n+1} \times \left( \mathbb{Z}^n \rtimes_{\varphi} \operatorname{GL}_n(\mathbb{Z}) \right)
\]
where $\varphi \colon \operatorname{GL}_n(\mathbb{Z}) \to \operatorname{Aut}(\mathbb{Z}^n)$ is the obvious isomorphism and the actions of
\begin{align*}
    \mathbf{v} &= (v_0,\dots,v_n) \in (R^{\times})^{n+1} \\
    \mathbf{w} &= (w_1,\dots,w_n) \in \mathbb{Z}^n \\
    M &= \begin{pmatrix} \mathbf{m}_1 \\ \vdots \\ \mathbf{m}_n \end{pmatrix}  \in \operatorname{GL}_n(\mathbb{Z})
\end{align*}
on a Laurent polynomial $P \in R[x_1^{\pm 1},\dots,x_n^{\pm 1}]$ are given by
\begin{align*}
    \mathbf{v} \ast P &:= v_0 \, P(v_1 \, x_1,\dots,v_n \, x_n) \\
    \mathbf{w} \ast P &:= \underline{x}^{\mathbf{w}} \, P \\
    M \ast P &:= P(\underline{x}^{\mathbf{m}_1},\dots,\underline{x}^{\mathbf{m}_n}).
\end{align*}
where for every $\mathbf{z} \in \mathbb{Z}^n$ we define $\underline{x}^{\mathbf{z}} := x_1^{z_1} \cdots x_n^{z_n}$. 
For every $\alpha \in (R^{\times})^{n+1} \times \operatorname{Aut}(\mathbb{G}_{m}^n)$ we will write $P_{\alpha} := \alpha \ast P$ and $\mathbf{v}_{\alpha} \in (R^{\times})^{n+1}$, $\mathbf{w}_{\alpha} \in \mathbb{Z}^n$, $M_{\alpha} \in \operatorname{GL}_n(\mathbb{Z})$ are the corresponding components.

\begin{remark} \label{rem:changing_functions}
    Let $\kappa$ be a field, let $P \in \kappa[x^{\pm 1},y^{\pm 1}]$ and let $\alpha \in (\kappa^{\times})^{3} \times \operatorname{Aut}(\mathbb{G}_{m}^2)$. Then we have an isomorphism $V_P \xrightarrow{\sim} V_{P_{\alpha}}$ between the zero loci of $P$ and $P_{\alpha}$ inside $\mathbb{G}_m^2$. This induces an isomorphism $\kappa(\widetilde{V_{P_{\alpha}}}) \xrightarrow{\sim} \kappa(\widetilde{V_P})$ between the function fields of the desingularisations of their compactifications, which identifies the functions $x,y \in \kappa(\widetilde{V_P})$ with $x_{\alpha} := x^a \, y^b$ and $y_{\alpha} := x^c \, y^d$, where $a,b,c,d \in \mathbb{Z}$ are such that $M_{\alpha} = \begin{psmallmatrix} a & b \\ c & d \end{psmallmatrix}$.
    
    Let now $J_P := J_{\widetilde{V_P}}$ denote the Jacobian of $\widetilde{V_P}$ (see \cref{not:diamond_operator}), let $G \leq J_P(\overline{\kappa})$ denote any subgroup such that $S_{x,y} \subseteq G$ and let $\psi \colon \mathbb{Q}[G] \to \mathbb{Q}$ be any $\mathbb{Q}$-linear map which is odd, \textit{i.e.} such that $\psi((-x)) = - \psi((x))$ for every $x \in G$. Then we have that $S_{x_{\alpha},y_{\alpha}} = S_{x,y}$ and
    \begin{align}
        \begin{pmatrix} \operatorname{div}(x_{\alpha}) \\ \operatorname{div}(y_{\alpha})
        \end{pmatrix} &= M_{\alpha} \,
        \begin{pmatrix} \operatorname{div}(x) \\ \operatorname{div}(y)
        \end{pmatrix}
        \label{eq:divisors_transformed_functions} \\
        \psi(\operatorname{div}(x_{\alpha}) \diamondsuit \operatorname{div}(y_{\alpha})) &= \det(M_{\alpha}) \, \psi(\operatorname{div}(x) \diamondsuit \operatorname{div}(y))
        \label{eq:rational_factor_transformed_functions}
    \end{align}
    which follows simply from the fact that $\diamondsuit$ is bilinear and that $\psi$ is odd.
\end{remark}

Before moving on, let us introduce some last pieces of notation.

\begin{notation}
    For every Laurent polynomial $P \in \mathbb{C}[x^{\pm 1},y^{\pm 1}]$ we denote by $V_P \hookrightarrow \mathbb{G}_m^2$ its zero locus, by $\overline{V_P} \hookleftarrow V_P$ a compactification of $V_P$ and by $\widetilde{V_P} \twoheadrightarrow \overline{V_P}$ a desingularization of $\overline{V_P}$. Moreover, we denote by $V_P^{\text{sing}} \hookrightarrow V_P$ the closed subset of singular points and by $V_P^{\text{reg}} := V_P \setminus V_P^{\text{sing}}$ its open complement. 
\end{notation}

\begin{definition}[see {\cite[Assumptions~3.2]{de97}}] \label{def:deninger_cycle}
    Let $P \in \mathbb{C}[x^{\pm 1},y^{\pm 1}]$ be any Laurent polynomial. Then we define a path
    \[
        \gamma_P := \{ (x,y) \in V_P(\mathbb{C}) \mid \lvert x \rvert = 1, \ \lvert y \rvert \leq 1 \}
    \]
    and we denote by $[\gamma_P] \in H_1(V_P(\mathbb{C}),\partial \gamma_P;\mathbb{Z})$ its class in singular homology.
\end{definition}

\begin{remark}[Amoeba map] \label{rmk:amoeba}
    The \textit{amoeba map} $\mu \colon \mathbb{G}_m^2(\mathbb{C}) \to \mathbb{R}^2$ is defined by
    \begin{align*}
        \mu \colon \mathbb{G}_m^2(\mathbb{C}) &\to \mathbb{R}^2 \\
        (x,y) &\mapsto (\log\lvert x \rvert, \log\lvert y \rvert)
    \end{align*}
    and it deserves this name because for every Laurent polynomial $Q \in \mathbb{C}[x^{\pm 1},y^{\pm 1}]$ the set $\mu(V_Q(\mathbb{C})) \subseteq \mathbb{R}^2$ is given by a bounded region to which are attached some ``tentacles'' going towards infinity (see \cite[Page~194]{gkz94} for a picture).
    In particular, the complement $\mathbb{R}^2 \setminus \mu(V_Q(\mathbb{C}))$ has at least one unbounded connected component (see \cite[Corollary~6.1.8]{gkz94}).
\end{remark}

Using the action of $(R^{\times})^{n+1} \times \operatorname{Aut}(\mathbb{G}_{m}^n)$ we can transform any Laurent polynomial to make the Deninger path (see \cref{def:deninger_cycle}) avoid the unit torus and the set of singular points. This can be done combining work of Besser and Deninger (see \cite[Fact~2.1]{bd99}) and Bornhorn (see \cite[Lemma~5.2.8]{bor99} and \cite[Lemma~1.7]{bor15}).

\begin{lemma} \label{lem:changing_the_polynomial}
    Let $Q \in \mathbb{Q}[x^{\pm 1},y^{\pm 1}]$ be any Laurent polynomial. Then there exists $\alpha \in (\mathbb{Q}^{\times})^{3} \times \operatorname{Aut}(\mathbb{G}_{m}^2)$ such that $Q_{\alpha} \in 1 + y \, \mathbb{Z}[x,y]$, $V_{Q_{\alpha}}(\mathbb{C}) \cap \mathbb{T}^2 = \emptyset$ and $\gamma_{Q_{\alpha}} \cap V_{Q_{\alpha}}^{\text{sing}}(\mathbb{C}) = \emptyset$, where $\mathbb{T}^2 \subseteq \mathbb{G}_m^2(\mathbb{C})$ denotes the real unit torus.
\end{lemma}
\begin{proof}
    First of all observe that we can write $Q = x^a y^b \, \left( (y^c (k + y \, Q_2)) + x \, Q_1 \right)$ for some $a,b,c \in \mathbb{Z}$, $k \in \mathbb{Q}^{\times}$ and $Q_1, Q_2 \in \mathbb{Q}[x,y]$. 
    Indeed, first of all we can write $Q = x^a y^b \, Q'$ for some $Q' \in \mathbb{Q}[x,y]$ and $a,b \in \mathbb{Z}$ such that $x, y \nmid Q'$.
    Hence we can write $Q' = Q'' + x \, Q_1$ for some $Q'', Q_1 \in \mathbb{Q}[x,y]$. 
    Finally, there exists $c \in \mathbb{Z}$ such that $Q'' = y^c \, Q'''$ for some $Q''' \in \mathbb{Q}[x,y]$, which implies that we can write $Q''' = k + y \, Q_2$ for some $k \in \mathbb{Q}^{\times}$ and $Q_2 \in \mathbb{Q}[x,y]$.
    
    Now, observe that there exist $p,q \in \mathbb{Q}^{\times}$ such that if
    $\tilde{Q}(x,y) := Q(p \, x,q \, y)$ then 
    $V_{\tilde{Q}}(\mathbb{C}) \cap \mathbb{T}^2 = \emptyset$. 
    To show this we can use the amoeba map $\mu \colon \mathbb{G}_m^2(\mathbb{C}) \to \mathbb{R}^2$ (see \cref{rmk:amoeba}). Indeed, $V_{\tilde{Q}}(\mathbb{C}) \cap \mathbb{T}^2 = \emptyset$ is equivalent to say that $0 \not\in \mu(V_{\tilde{Q}}(\mathbb{C}))$. Moreover, we know that $\mu(V_{\tilde{Q}}(\mathbb{C})) = \tau_{p,q}(\mu(V_Q(\mathbb{C})))$, where $\tau_{p,q} \colon \mathbb{R}^2 \to \mathbb{R}^2$ denotes the translation by the vector $-(\log\lvert p \rvert,\log\lvert q \rvert)$.
    Hence we can use the fact that $\mathbb{R}^2 \setminus \mu(V_Q(\mathbb{C}))$ has at least one unbounded connected component to see that there exist $p,q \in \mathbb{Q}^{\times}$ sufficiently large such that $V_{\tilde{Q}}(\mathbb{C}) \cap \mathbb{T}^2 = \emptyset$.
    
    Now, let us write $\mu(V_Q^{\text{sing}}(\mathbb{C})) = \{ (\lambda_1^{(i)}, \lambda_2^{(i)}) \}_{i = 1}^s$ and $Q = 1 + \sum_{j = 1}^n a_j(x) \, y$ for some polynomials $\{a_j\}_{j = 1}^n \subseteq \mathbb{Z}[x]$. Let finally $m \in \mathbb{N}$ be any natural number such that
    \[
        m \geq \max\left( \{ \deg_x(a_j) \}_{j = 1}^n \cup \left\{ \left\lvert \lambda_2^{(i)}/\lambda_1^{(i)} \right\rvert \, \colon \, \lambda_1^{(i)} \neq 0 \right\}_{i = 1}^s \right)
    \]
    and let us take $\alpha \in (\mathbb{Q}^{\times})^3 \times \operatorname{Aut}(\mathbb{G}_m^2)$ to be
    \[
        \mathbf{v}_{\alpha} = \left( \frac{1}{k}, p,q \right), \quad \mathbf{w}_{\alpha} = (-a, -b-c) \quad \text{and} \quad M_{\alpha'} = \begin{pmatrix} 1 & c+1 \\ 0 & 1 \end{pmatrix} \, \begin{pmatrix} -1 & -1 \\ m & m+1 \end{pmatrix} .
    \]
    
    Then the fact that $Q_{\alpha} \in 1 + y \, \mathbb{Z}[x,y]$ and $V_{Q_{\alpha}}(\mathbb{C}) \cap \mathbb{T}^2 = \emptyset$ follow from the previous discussions, whereas $\gamma_{Q_{\alpha}} \cap V_{Q_{\alpha}}^{\text{sing}}(\mathbb{C}) = \emptyset$ follows from \cite[Lemma~5.2.8]{bor99} (see also \cite[Lemma~1.7]{bor15}).
\end{proof}

\begin{remark}
    If we start from a tempered polynomial $P \in \mathbb{Q}[x^{\pm 1},y^{\pm 1}]$ the resulting polynomial $P_{\alpha}$ will in general not be tempered anymore, because we are scaling its variables and therefore its coefficients. 
    Nevertheless, the functions $x_{\alpha}, y_{\alpha}$ will still be supported on torsion points, thanks to \eqref{eq:divisors_transformed_functions}, and $r_{\widetilde{V_P}}(\eta_{x_{\alpha},y_{\alpha}}) \neq 0$, thanks to \eqref{eq:rational_factor_transformed_functions}. Hence we will still be able to apply \cref{thm:rohrlich_general_theorem}, and we will find a relation between the Mahler measure of $P_{\alpha}$ and the $L$-value $L'(E,0)$ despite the fact that $P$ is not tempered.
\end{remark}

We are now ready to prove our main theorem.
\begin{theorem}[see \cref{thm:main_theorem}] \label{thm:main_theorem_body}
    Let $E$ be an elliptic curve defined over $\mathbb{Q}$ having potential complex multiplication by the ring of integers $\mathcal{O}_K$ of an imaginary quadratic field $K$.
    Let $f,g \in \mathbb{Q}(E)$ be any pair of functions such that $\mathbb{Q}(E) = \mathbb{Q}(f,g)$, $S_{f,g} \subseteq E(\overline{\mathbb{Q}})_{\text{tors}}$ and $\langle r_E^{\infty}(\eta_{f,g}), \gamma_E \rangle \neq 0$.
    Let $Q \in \mathbb{Q}[x^{\pm 1},y^{\pm 1}]$ be a minimal polynomial for $f,g$ and let $P := Q_{\alpha}$ for any $\alpha \in (\mathbb{Q}^{\times})^{3} \times \operatorname{Aut}(\mathbb{G}_{m}^2)$ satisfying the conditions of \cref{lem:changing_the_polynomial}.
    Let $c_P \in \mathbb{Q}^{\times}$ be defined by the identity $[\widetilde{j}_{\ast}(j_{\ast}(\gamma_P))] = (- c_P) \, \gamma_E$, where $j$ denotes the open embedding $j \colon V_P^{\text{reg}} \hookrightarrow E \setminus S_{f,g}$ and $\widetilde{j}$ denotes the open embedding $\widetilde{j} \colon E \setminus S_{f,g} \hookrightarrow E$.
    Let moreover $\{[\gamma_z]\}_{z \in S_{f,g}} \subseteq H_1(E(\mathbb{C}) \setminus S_{f,g};\mathbb{Z})$ be the homology classes associated to small loops around each point $z \in S_{f,g}$, let $[\beta_1], [\beta_2] \in H_1(E(\mathbb{C}) \setminus S_{f,g};\mathbb{Z})$ be such that the set $\{[\beta_i]\}_{i = 1}^2 \cup \{[\gamma_z]\}_{z \in S_{f,g}}$ generates $H_1(E(\mathbb{C}) \setminus S_{f,g};\mathbb{Z})$ and let $\{a_z\}_{z \in S_{f,g} \setminus \{0\}} \subseteq \mathbb{Z}$ and $\{b_i\}_{i = 1}^2 \subseteq \mathbb{Z}$ be defined by the decomposition
    \[
            [j_{\ast}(\gamma_{P})] = \sum_{z \in S_{f,g} \setminus \{0\}} a_z \, [\gamma_z] + \sum_{i = 1}^2 b_i \, [\beta_i]
    \]
    which exists and is unique thanks to the exact sequence \eqref{eq:homology_punctured_curve}.
    Finally, define
    \begin{equation} \label{eq:explicit_numbers}
        \begin{aligned}
            r &:= \frac{c_P \, \mathcal{R}(\operatorname{div}(\tilde{x}) \diamondsuit \operatorname{div}(\tilde{y}))}{n_{\tilde{x},\tilde{y}}} \in \mathbb{Q} \\
            s &:= \prod_{z \in S_{\tilde{x},\tilde{y}} \setminus \{ 0 \}} \partial_z(\{\tilde{x},\tilde{y}\})^{a_z} \in \overline{\mathbb{Q}}^{\times}
        \end{aligned}
    \end{equation}
    where $\tilde{x}, \tilde{y} \in \mathbb{Q}(E)$ are given by $\tilde{x} := f_{\alpha}$ and $\tilde{y} := g_{\alpha}$ (see \cref{rem:changing_functions}).
    Then $V_P$ is birational to $E$ and
    \[
    m(P) = r \, L'(E,0) + \log\lvert s \rvert \]
    with $r \neq 0$ for a suitable choice of $\alpha$.
\end{theorem}
\begin{proof}
    Recall first of all that $S_{\tilde{x},\tilde{y}} = S_{f,g}$.
    Observe moreover that 
    \[
    [\widetilde{j}_{\ast}(j_{\ast}(\gamma_P))],\gamma_E \in H_1(E(\mathbb{C});\mathbb{Q})^{-} \cong \mathbb{Q}
    \]
    which implies that $c_P \in \mathbb{Q}$ exists.
    We have now the following chain of identities:
    \begin{align}
        \label{eq:main_theorem_1}
        m(P) &= - \, \left\langle r_{V_{P}}^{\infty}(\{x,y\}), [\gamma_{P}] \right\rangle = \\
        \label{eq:main_theorem_2}
        &= - \, \left\langle r_{V_{P}^{\text{reg}}}^{\infty}(\iota^{\ast}(\{x,y\})), [\gamma_{P}] \right\rangle = \\
        \label{eq:main_theorem_3}
        &= - \, \left\langle r_{E \setminus S_{\tilde{x},\tilde{y}}}^{\infty}(\{\tilde{x},\tilde{y}\}), [j_{\ast}(\gamma_{P})] \right\rangle = \\
        \label{eq:main_theorem_4}
        &= - \, \frac{1}{n_{\tilde{x},\tilde{y}}} \, \left( \left\langle r_E^{\infty}(\eta_{\tilde{x},\tilde{y}}),[\widetilde{j}_{\ast}(j_{\ast}(\gamma_{P}))] \right\rangle - \sum_{z \in S_{\tilde{x},\tilde{y}} \setminus \{ 0 \}} \left\langle r_{E \setminus S_{\tilde{x},\tilde{y}}}^{\infty}(\{ \partial_z(\{\tilde{x},\tilde{y}\}),\varphi_{\tilde{x},\tilde{y}}^{(z)} \}),[j_{\ast}(\gamma_{P})] \right\rangle  \right) = \\
        \label{eq:main_theorem_5}
        &= \frac{c_{P}}{n_{\tilde{x},\tilde{y}}} \, \left\langle r_E^{\infty}(\eta_{\tilde{x},\tilde{y}}),\gamma_E \right\rangle + \sum_{z \in S_{\tilde{x},\tilde{y}} \setminus \{ 0 \}} a_z \, \log\lvert \partial_z(\{\tilde{x},\tilde{y}\}) \rvert = \\ 
        \label{eq:main_theorem_6}
        &= \left( \frac{c_P \, \mathcal{R}(\operatorname{div}(\tilde{x}) \diamondsuit \operatorname{div}(\tilde{y}))}{n_{\tilde{x},\tilde{y}}} \right) \, L'(E,0) + \log\left\lvert \prod_{z \in S_{\tilde{x},\tilde{y}} \setminus \{ 0 \}} \partial_z(\{\tilde{x},\tilde{y}\})^{a_z} \right\rvert
    \end{align}
    where $\iota$ denotes the open embedding $\iota \colon V_{P}^{\text{reg}} \hookrightarrow V_{P}$.
    
    To explain these identities we observe that \eqref{eq:main_theorem_1} is an application of \cite[Theorem~3.4]{de97}, using the fact that $P^{\ast} = 1$ and \eqref{eq:main_theorem_2} is a consequence of the fact that $\gamma_{P} \subseteq V_{P}^{\text{reg}}(\mathbb{C})$. Moreover, \eqref{eq:main_theorem_3} follows from the fact that $\iota^{\ast}(\{x,y\}) = j^{\ast}(\{\tilde{x},\tilde{y}\})$ and \eqref{eq:main_theorem_4} follows from the definition of $\eta_{\tilde{x},\tilde{y}}$. Finally, \eqref{eq:main_theorem_5} follows from \cref{prop:regulator_open_curve} and \eqref{eq:main_theorem_6} follows from \cref{thm:rohrlich_general_theorem}.
    
    Now observe that $c_P = 0$ if $V_P(\mathbb{C}) \cap \{(x,y) \in \mathbb{C}^2 \colon \lvert x \rvert = 1 \} \subseteq \{(x,y) \in \mathbb{C}^2 \colon \lvert y \rvert < 1 \}$ (see \cite[Page~48]{bo98}). 
    Clearly, the same holds if we take $\lvert y \rvert > 1$ in the set on the right and if we change $x$ with $y$. In other words, if the amoeba $\mu(V_P(\mathbb{C}))$ does not intersect all the four semi-axes we have that $c_P = 0$.
    Nevertheless, it is clear that we can translate the amoeba sufficiently enough so that, with a convenient rotation, it will intersect all the four semi-axes. When this happens, we will have that $c_P \neq 0$.
\end{proof}
\begin{remark}
    Pairs of functions like the ones described in the statement of \cref{thm:main_theorem_body} are given by the constructions of Deninger and Wingberg (see \cref{lem:deninger_wingberg_construction}) and Rohrlich (see \cref{lem:rohrlich_functions}).
\end{remark}

\begin{appendices}

\crefalias{section}{appsec}
\titleformat{\section}[frame]
{\normalfont}
{\filright
\footnotesize
\enspace APPENDIX \thesection\enspace}
{8pt}
{\Large\bfseries\filcenter}

\section{Conductors of abelian varieties with complex multiplication}
\label{sec:conductors}

The aim of this appendix is to provide references for the theory of complex multiplication, and for a proof of the following result, which is due to Deuring.

\begin{proposition}[Deuring] \label{prop:deuring_CM_conductor}
    Let $E$ be an elliptic curve defined over $\mathbb{Q}(j)$ where $j = j(E)$. Suppose that $\operatorname{End}_{\overline{\mathbb{Q}}}(E) \cong \mathcal{O}$ for some order $\mathcal{O} \subseteq \mathcal{O}_F$ inside an imaginary quadratic field $F$. Then we have that 
    \[
        \mathfrak{f}_E = \operatorname{N}_{F(j)/\mathbb{Q}(j)}(\mathfrak{f}_{\psi_E}) \,  \operatorname{disc}(F(j)/\mathbb{Q}(j))
    \]
    where $\mathfrak{f}_E \subseteq \mathcal{O}_{\mathbb{Q}(j)}$ denotes the conductor ideal of $E$, $\psi_E \colon \mathbb{A}_{F(j)}^{\times} \to \mathbb{C}^{\times}$ denotes the Hecke character associated to $E$, with conductor $\mathfrak{f}_{\psi_E} \subseteq \mathcal{O}_{F(j)}$, and $\operatorname{disc}(F(j)/\mathbb{Q}(j)) \subseteq \mathcal{O}_{\mathbb{Q}(j)}$ denotes the discriminant of the extension $\mathbb{Q}(j) \subseteq F(j)$.
\end{proposition}

Let us recall first of all the notion of complex multiplication (see \cite[Chapter~1]{cco14} for an excellent introduction). 
Every abelian variety $A$ defined over a field $\kappa$ is isogenous to a unique product $A_1^{e_1} \times \dots \times A_n^{e_n}$, where each $A_i$ is \textit{simple}, \textit{i.e.} it does not have any non-trivial abelian sub-variety.
Hence we have an isomorphism of $\mathbb{Q}$-algebras 
\[
    \operatorname{End}(A)_{\mathbb{Q}} := \operatorname{End}(A) \otimes_{\mathbb{Z}} \mathbb{Q} \cong \prod_{i = 1}^n \operatorname{Mat}_{e_i}(\operatorname{End}(A_i)_{\mathbb{Q}})
\]
and every choice of polarisation $A_i \to A_i^{\vee}$ endows the $\mathbb{Q}$-algebra $\operatorname{End}(A_i)_{\mathbb{Q}}$ with a positive involution. 
Since $A_i$ is simple then $\operatorname{End}(A_i)_{\mathbb{Q}}$ is a \textit{simple division algebra}, \textit{i.e.} $\operatorname{End}(A_i)_{\mathbb{Q}}$ does not have any non-trivial two-sided ideal and for every $\alpha, \beta \in \operatorname{End}(A_i)_{\mathbb{Q}}$ there exists a unique pair $\gamma, \delta \in \operatorname{End}(A_i)_{\mathbb{Q}}$ such that $\alpha = \gamma \, \beta = \beta \, \delta$.
Hence the algebras $\operatorname{End}(A_i)_{\mathbb{Q}}$ fall within Albert's classification of division $\mathbb{Q}$-algebras with a positive involution (see \cite[\S~1.3.6]{cco14}).

Recall now that if $M$ is a simple algebra over a field $\kappa$ then its center $\operatorname{Z}(M)$ is isomorphic to a field $\kappa' \supseteq \kappa$ and $[M \colon Z(M)]$ is a square. Hence we can define the \textit{reduced degree} 
\[
    [M \colon \kappa]^{\text{red}} := \sqrt{[M \colon Z(M)]} \, [\kappa' \colon \kappa] \in \mathbb{N}
\]
and if $N = \prod_i N_i$ is a semi-simple $\kappa$-algebra we define $[N \colon \kappa]^{\text{red}} := \sum_i [N_i \colon \kappa]^{\text{red}}$.
Using this notation, we have that $[\operatorname{End}(A)_{\mathbb{Q}} \colon \mathbb{Q}]^{\text{red}} \leq 2 \dim(A)$ and we know that $[\operatorname{End}(A)_{\mathbb{Q}} \colon \mathbb{Q}]^{\text{red}} = 2 \dim(A)$ if and only if for every $i \in \{ 1, \dots, n \}$ there exists a totally imaginary number field $F_i$ which contains a totally real sub-field $F_i^{+} \subseteq F_i$ such that $[F_i \colon F_i^{+}] = 2$ and $F_i \hookrightarrow \operatorname{End}(A_i)_{\mathbb{Q}}$. 
Such types of number fields are called \textit{CM fields}, in view of the following definition.
\begin{definition} \label{def:complex_multiplication}
    Let $A$ be an abelian variety over a field $\kappa$. Then we say that $A$ has \textit{complex multiplication} if $[\operatorname{End}(A)_{\mathbb{Q}} \colon \mathbb{Q}]^{\text{red}} = 2 \dim(A)$ and we say that it has \textit{potential complex multiplication} if there exists a finite extension $\kappa' \supseteq \kappa$ such that $A_{\kappa'}$ has complex multiplication.
\end{definition}
\begin{remark}
    Let $A$ be an abelian variety over a field $\kappa$, which has complex multiplication.
    Then we know from Albert's classification that either $\operatorname{End}(A_i)_{\mathbb{Q}} = F_i$ or $\operatorname{char}(\kappa) > 0$ and $\operatorname{End}(A_i)_{\mathbb{Q}}$ is a non-split quaternion algebra.
\end{remark}

We see that if an abelian variety $A$ defined over a field $\kappa$ has potential complex multiplication then we have an embedding $E \hookrightarrow \operatorname{End}(A_{\overline{\kappa}})_{\mathbb{Q}}$, where $E$ is a \textit{CM algebra}, \textit{i.e.} a product of CM fields. 
If $\operatorname{char}(\kappa) = 0$ we have an action of $E$ on the tangent space of $A_{\overline{\kappa}}$ at the origin. 
This determines a \textit{CM type} of $E$, \textit{i.e.} a collection of algebra homomorphisms $\Phi \subseteq \operatorname{Hom}_{\mathbb{Q}}(E,\mathbb{C})$ such that $\Phi \cap \overline{\Phi} = \emptyset$ and $\Phi \cup \overline{\Phi} = \operatorname{Hom}_{\mathbb{Q}}(E,\mathbb{C})$, where $\overline{\Phi}$ is obtained from $\Phi$ by composing with complex conjugation.
If $E$ is a CM algebra and $\Phi \subseteq \operatorname{Hom}_{\mathbb{Q}}(E,\mathbb{C})$ is a CM type we call $(E,\Phi)$ a \textit{CM pair}.
Moreover, $E$ acts as well on the spaces of differential forms defined on $A$, and this action can be used to study the field where $A$ acquires complex multiplication. 
This is summarised in the following proposition.
 
\begin{proposition}
    Let $A$ be an abelian variety defined over a field $K$ of characteristic zero, which has complex multiplication by a CM algebra $E = F_1 \times \dots \times F_r$. Let $\Phi_i$ be the CM type induced on each $F_i$, and let $(F_i^{\ast},\Phi_i^{\ast})$ be the reflex CM pairs (see  \cite[Section~8.3]{sh98}). 
    Then there is an embedding $F_i^{\ast} \hookrightarrow K$ for every $i \in \{ 1, \dots, r\}$.
    
    Conversely, suppose that $A$ is a simple abelian variety defined over a field $K$ of characteristic zero, such that $\operatorname{End}(A_{K'}) \cong F$ for some finite extension $K' \supseteq K$ and some CM field $F$. Then $\operatorname{End}(A_{L})_{\mathbb{Q}} = \operatorname{End}(A_{K'})_{\mathbb{Q}}$ for every sub-field $L \subseteq K'$ such that $L \supseteq K$ and $L \supseteq \iota(F^{\ast})$. Here $(F^{\ast},\Phi^{\ast})$ is the reflex CM pair of $(F,\Phi)$, where $\Phi$ is the CM type induced by complex multiplication on $A_{K'}$, and $\iota \colon F \hookrightarrow K'$ is the embedding given in the previous paragraph.
\end{proposition}
\begin{proof}
    See \cite[Chapter~II,~Proposition~30]{sh98}.
\end{proof}

We can now come to the issue of relating the conductor of an abelian variety with potential complex multiplication to the conductor of the corresponding Hecke character. 
The main theorem that we are going to use is the following one, which is essentially due to Milne.

\begin{theorem}[Milne]
    Let $A$ be an abelian variety defined over a number field $K$, let $L \supseteq K$ be a finite Galois extension and suppose that $A_L$ has complex multiplication by the CM algebra $E \hookrightarrow \operatorname{End}(A_L)_{\mathbb{Q}}$. 
    Let $d = \dim(A)$, $m = [L \colon K]$ and assume that $E \cap \operatorname{End}(A)_{\mathbb{Q}}$ is a field and that $[E \colon E \cap \operatorname{End}(A)_{\mathbb{Q}}] = m$.
    
    Then we have that $m \mid 2 d$ and for every $\sigma \colon L \hookrightarrow \mathbb{C}$ we have that $\mathfrak{f}_A = \left( \operatorname{N}_{L/K}(\mathfrak{f}_{\chi_{\sigma}}) \, \operatorname{disc}(L/K) \right)^{2 d/m}$, where $\chi_{\sigma}$ denotes the Hecke character associated to $A_L$ and $\sigma$ (see \cite[Section~7]{st68}), $\mathfrak{f}_A$ denotes the conductor of $A$ (see \cite[\S~1.(b)]{mi72}) and $\mathfrak{f}_{\chi_{\sigma}}$ denotes the conductor of the Hecke character $\chi_{\sigma}$ (see \cite[Section~VII.6]{ne99}). 
\end{theorem}
\begin{proof}
    Since $[E \colon \mathbb{Q}] = 2d$ it is immediate to see that $m \mid 2d$. Then the theorem follows from \cite[Theorem~3]{mi72} and from the two formulas $\mathfrak{f}_{A_L} = \mathfrak{f}_{\chi_{\sigma}}^{2 d}$ (see \cite[Theorem~12]{st68}) and $\mathfrak{f}_{\operatorname{N}_{L/K}(B)} = \operatorname{N}_{L/K}(\mathfrak{f}_B) \, \operatorname{disc}(L/K)^{2 d}$ (see \cite[Theorem~1]{mi72}). 
    In the second formula, $\operatorname{N}_{L/K}(B)$ denotes the Weil restriction of an abelian variety $B$ defined over $L$ (see \cite[Section~7.6]{blr90} and \cite[Section~A.5]{cgp15}).
\end{proof}

As Milne already states in \cite{mi72}, this theorem applies in particular when $A$ is simple (over $K$) and $L$ is the smallest Galois extension of $K$ such that $\operatorname{End}(A_L)_{\mathbb{Q}}$ contains the center of $\operatorname{End}(A_{\overline{K}})_{\mathbb{Q}}$. This is exactly the situation of \cref{prop:deuring_CM_conductor}, which gives us a modern proof of Deuring's result.
\end{appendices}

\section*{Acknowledgements}

We would like to thank our advisers Ian Kiming and Fabien Pazuki for their unceasing support and for the many mathematical discussions around this project. 
We would like to thank Fran\c{c}ois Brunault and Wadim Zudilin for having invited us to visit them, and for many useful discussions on this project. 
We would like to thank our colleague Francesco Campagna for the many clarifications of small and big doubts given in the course of this project.
Finally, we would like to thank Lars Hesselholt and Michalis Neururer for useful discussions.

This work is part of the author's PhD thesis at the University of Copenhagen, which we would like to thank for their financial support. Moreover, we would like to acknowledge the hospitality of the {\'E}cole normale sup{\'e}rieure de Lyon and the Radboud University in Nijmegen, where some of this research was conducted.

\printbibliography

\end{document}